\documentclass{article}

\usepackage{a4, amsthm, amssymb, amscd, mathptmx}

\DeclareSymbolFont{rsfs}{U}{rsfs}{m}{n}
\DeclareSymbolFontAlphabet{\mathscr}{rsfs}

\newcommand{\catA}{\ul{\hbox{\rm A}}}
\newcommand{\catg}{\ul{\hbox{\rm g}}}
\newcommand{\catM}{\ul{\hbox{\rm M}}}
\newcommand{\catP}{\ul{\hbox{\rm P}}}
\newcommand{\catQ}{\ul{\hbox{\rm Q}}}
\newcommand{\catR}{\ul{\hbox{\rm R}}}
\newcommand{\catT}{\ul{\hbox{\rm T}}}
\newcommand{\catt}{\ul{\hbox{\rm t}}}

\newcommand{\Coh}{\ul{\mathop{\rm Coh}}}
\newcommand{\Rep}{\ul{\mathop{\rm Rep}}}
\newcommand{\catq}{\ul{\mathop{\rm Q}}}
\newcommand{\catp}{\ul{\mathop{\rm P}}}
\newcommand{\cata}{\ul{\mathop{\rm A}}}

\renewcommand{\Bbb}{\mathbb}
\renewcommand{\frak}{\mathfrak}
\renewcommand{\cal}{\mathscr}
\newcommand{\catqot}{/\hskip-3pt/}
\newcommand{\C}{{\Bbb C}}

\newcommand{\E}{{\cal E}}
\newcommand{\ED}{\mathop{\rm END}}
\newcommand{\End}{\mathop{\rm End}}

\newcommand{\F}{{\cal F}}
\newcommand{\G}{{\cal G}}

\newcommand{\GL}{\mathop{\rm GL}}

\newcommand{\gr}{\mathop{\rm gr}}

\newcommand{\Hom}{\mathop{\rm Hom}}

\renewcommand{\Im}{\mathop{\rm Im}}
\newcommand{\In}{\mathop{\rm In}}
\newcommand{\id}{\mathop{\rm id}}
\renewcommand{\L}{{\cal L}}

\newcommand{\n}{{\cal N}}
\renewcommand{\O}{{\cal O}}
\newcommand{\Out}{\mathop{\rm Out}}
\renewcommand{\P}{{\Bbb P}}
\newcommand{\Pic}{\mathop{\rm Pic}}

\newcommand{\pt}{\mathop{\rm pt}}
\newcommand{\Q}{{\Bbb Q}}

\newcommand{\SL}{\mathop{\rm SL}}
\newcommand{\Spec}{\mathop{\rm Spec}}

\renewcommand{\tilde}{\widetilde}
\newcommand{\Z}{{\Bbb Z}}
\newcommand{\N}{{\Bbb N}}
\newcommand{\la}{\lambda}
\newcommand{\lra}{\longrightarrow}
\newcommand{\ra}{\rightarrow}

\newcommand{\lma}{\longmapsto}
\newcommand{\p}{\prime}
\newcommand{\q}{\quad}

\renewcommand{\phi}{\varphi}
\newcommand{\rk}{\mathop{\rm rk}}
\newcommand{\eps}{\varepsilon}
\newcommand{\St}{\mathop{\rm Star}}

\newcommand{\ul}{\underline}
\newcommand{\ol}{\overline}

\theoremstyle{plain}
\newtheorem{Thm}{Theorem}[section]
\newtheorem{Cor}[Thm]{Corollary}
\newtheorem{Prop}[Thm]{Proposition}
\newtheorem{Lem}[Thm]{Lemma}

\newtheorem*{Thm*}{Theorem}

\theoremstyle{remark}
\newtheorem{Rem}[Thm]{Remark}
\newtheorem{Ass}[Thm]{Assumptions}
\newtheorem*{Conv}{Convention}
\newtheorem{Ex}[Thm]{\it Example}

\begin{document}


\title{\sf Moduli problems of sheaves associated with oriented trees}
\author{Alexander Schmitt}
\maketitle 
\begin{abstract}
To every oriented tree we associate
vector bundle problems. 
We define semistability concepts for these
vector bundle problems and establish the existence of moduli spaces.
As an important application, we obtain an algebraic construction
of the moduli space of holomorphic triples.
\end{abstract}

\section*{Introduction}

Let  $Q=(V,A,t,h)$ be an oriented graph, or quiver, 
and $X$ a fixed projective algebraic manifold. Associate
to each vertex $i\in V$ a coherent sheaf $\E_i$ over $X$ and to each
arrow a homomorphism $\phi_a\in \Hom(\E_{t(a)}, \E_{h(a)})$, and let
all the sheaves $\E_i$, $i\in V$, and homomorphisms $\phi_a$, $a\in A$,
vary. This yields interesting moduli problems in Algebraic Geometry.
\par
To the author's knowledge,
this kind of problems has been studied only in a few special cases: 
First, the results of King \cite{Ki} cover the case when $X$ is a point,
second,
the theory of semistable torsion free coherent sheaves deals with  
the quiver $\bullet$ (cf.\ \cite{HL}).
Apart from this, the objects associated with 
$\bullet\lra\bullet$ when $X$ is a curve, the so-called
holomorphic triples (which we will review below)
were treated by Bradlow and Garcia-Prada
\cite{BG}, \cite{GP}, the objects associated with $A_n$-quivers
($\bullet\lra\cdots\lra \bullet$) have been studied under the name
of holomorphic chains in the context of dimensional reduction
by \'Alvarez-C\'onsul and Grac\'\i a-Prada \cite{AG}, and, 
finally, Higgs bundles
or Hitchin pairs \cite{Hi}, \cite{Ni}, \cite{Si}, \cite{Yo}, \cite{Sch}
can be seen as
the objects associated with the quiver consisting of one vertex joined
by an arrow to itself.
More recently, \'Alvarez-C\'onsul and Grac\'\i a-Prada have announced
generalizations of their work \cite{AG} to quivers without oriented
cycles.
\par
In this note, we address the case when $X$ is an
arbitrary projective manifold and $Q$ is an oriented tree.
One of the central points of the paper is the definition of  
$\vartheta$-semistability
for representations of $Q$, i.e., tuples $(\E_i,i\in V;\phi_a, a\in A)$
as above.
The notion we propose is obtained by applying principles from our
paper \cite{Sch3}. It also fits in the framework suggested by King \cite{Ki}
for finding the semistability concept.
The main result of our paper is the existence of moduli spaces for
the $\vartheta$-semistable representations.
The construction of the moduli spaces itself applies the
general GIT machinery of \cite{GIT}, adapted to vector bundle
problems by Seshadri, Gieseker, Maruyama, Simpson and others
(see \cite{HL} and \cite{Sch3} for precise references).
Of course, due to the complexity of the objects we study, the
details become technically quite involved.
\par
The motivation to study these questions is the following:
First, the case of $X=\{\pt\}$, i.e., King's work, is
important because of its relation to the classification
theory of algebras \cite{Ga}, \cite{Ri}.
Recently, King's moduli spaces have found applications 
to non-commutative algebraic geometry
\cite{Ln}, and
generalizations of them 
were applied to the theory of quantum algebras
\cite{Lu}, \cite{Na1}, \cite{Na2}.
One might therefore hope that the moduli spaces we construct will find
applications in one or the other of these areas.
\par
Second, there is the special case of $\bullet\lra \bullet$ or, more
generally, the $A_n$-quivers.
Suppose $X$ is a smooth, projective curve. A holomorphic triple
is then a triple $(E_1,E_2,\phi)$ consisting of two vector bundles
$E_1$ and $E_2$ on $X$ and a homomorphism $\phi$ between them.
These objects were obtained by Bradlow and Garcia-Prada in \cite{BG}
and \cite{GP}
by a process of dimensional
reduction from certain $\mathop{\rm SU}(2)$-equivariant bundles on
$X\times\P_1$.
They defined the notion of $\tau$-(semi)stability
for holomorphic triples and gave a construction of the moduli space
of $\tau$-stable triples.
In \cite{BG}, it was suggested to
construct this moduli space via Geometric Invariant Theory.
Of course, the notion
of a holomorphic triple has an obvious interpretation on higher
dimensional manifolds, and our results contain the construction
on an arbitrary projective manifold $X$.
The central technical point in our GIT construction
is the identification
of the semistable points in a suitable
parameter space which is trickier than usual due to the failure
of a certain addivity property of the weights for the actions
one has to study. To solve this, we will prove a decomposition
theorem for one parameter subgroups of $\SL(V)\times\SL(W)$
w.r.t.\ the action on $\Hom(V,W)$.
It might also be interesting to note that the GIT construction
reveals that there are actually two parameters involved in the
definition of semistability for holomorphic triples.  
This might become useful for relating the moduli space of holomorphic triples
to other moduli spaces.
\par
Finally, the results of this paper are after the results of
the paper \cite{Sch3} a further step towards 
a universal theory working for (almost) all vector bundle
problems.
In such a theory, the input would be a representation of a
(reductive) algebraic group $G$
on a finite dimensional vector space $W$ --- which defines in a natural
way a moduli problem ---, and the output would be the
(parameter dependent) semistability concept and moduli spaces
for the semistable objects.
For a precise formulation in the case $G=\GL(r)$ and a solution
of the moduli problem over curves, we refer the reader
to our paper \cite{Sch3}.
The present paper corresponds to the case of $G=\prod_{i\in V}
\GL(r_i)$ with its representation on
$\bigoplus_{a\in A} \Hom(\C^{r_{t(a)}},\C^{r_{h(a)}})$.
\par
We work only with oriented trees, 
because they are well suited for inductions.
Indeed, they can be thought of as being built ``inductively'', just like 
a tree in the real world grows from a small tree to a big tree, spreading
more and more branches. Therefore, many of the problems reduce
to the quiver $\bullet\lra \bullet$.
Notably, an inductive procedure can be used to extend the already mentioned
decomposition result which is the key to the general theory.
At the moment, our techniques do not extend to arbitrary quivers, and examples
suggest that indeed further technical difficulties will arise for
quivers other than oriented trees.
\subsection*{Notations and conventions}
The base field is the field of complex numbers.
(The restriction to characteristic zero is necessary to apply
Maruyama's boundedness result (\cite{HL}, Thm.~3.3.7) for certain
families of torsion free coherent sheaves.)
\par
In the sequel, $X$ is understood to be a projective manifold,
and $\O_X(1)$ to be an ample invertible sheaf on $X$.
For any coherent sheaf $\E$, we write $P(\E)$ for
its Hilbert polynomial w.r.t.\ $\O_X(1)$, and
degrees and slopes are computed w.r.t.\ $\O_X(1)$.
Most of the time, we will fix some polynomials $P_i$, $i$ in some
index set.
Having done this, $r_i$, $d_i$, and $\mu_i$ are used for
the rank, degree, and slope defined by these polynomials
when interpreted as Hilbert polynomials w.r.t.\ $\O_X(1)$.
Recall that any torsion free coherent sheaf $\G$ has a uniquely
determined slope Harder-Narasimhan filtration
$0=\G_0\subset \G_1\subset\cdots\subset\G_l=\G$, such that all the
quotients $\G_i/\G_{i-1}$, $i=1,...,l$, are slope semistable,
and $\mu(\G_{i}/\G_{i-1})>\mu(\G_{i+1}/\G_i)$, $i=1,...,l-1$ 
(see \cite{HL}, \S 1.6). 
One sets $\mu_{\max}(\G):=\mu(\G_1)$ and $\mu_{\min}(\G):=
\mu(\G/\G_{l-1})$
Finally, $\chi(\G)=P(\G)(0)$ stands for the Euler
characteristic of $\G$.
\par
We use Grothendieck's convention for projectivizing a vector bundle $E$, 
i.e.,
$\P(E)$ is the bundle of hyperplanes in the fibres of the vector bundle $E$.
\subsection*{Acknowledgments}
The paper was written during the author's stays at the University
of Barcelona and at the Bar-Ilan University.
During this time, the author was supported by grant \#1996SGR00060
of the Generalitat de Catalunya and by a grant of the Emmy-Noether
Insitute, respectively.
\section{Statement of the main results}
\subsection*{Representations of quivers}
We first recall the notion of a representation of a quiver in an abelian
category and describe an abstract semistability concept, going back to King
\cite{Ki}.
\par
A \it quiver $Q$ \rm consists
of a set $V:=\{\, 1,...,n\,\}$ whose elements are called \it
vertices \rm and a set $A$, the \it arrows\rm,
together with maps $h,t\colon A\lra V$.
Given an arrow $a\in A$, we call $h(a)$ its \it head \rm 
and $t(a)$ its \it tail\rm.
\begin{Conv}
We will only deal with quivers without multiple arrows.
Therefore, we often write an arrow in the form $(t(a),h(a))$.
\end{Conv}
\begin{Ex} 
Given a quiver $Q$ as above, the \it underlying graph \rm
is the graph $\Gamma_Q$ whose set of vertices is just $V$
and whose set of edges is $E:=\bigl\{\, \{i_1,i_2\}\, |\,\allowbreak
(i_1,i_2)\in A\,\bigr\}$.
An \it oriented tree \rm is a quiver
$Q$ with connected underlying graph and
$\# V=\# A+1$. This means precisely that the underlying
graph $\Gamma_Q$ is a tree. 
\end{Ex}
A \it path in $Q$ from $i$ to $j$
\rm is a sequence of arrows $a_1,...,a_m$ with $t(a_1)=i$,
$h(a_l)=t(a_{l+1})$, $l=1,...,m-1$, and $h(a_m)=j$. Here, $m$ is called
\it the length of the path\rm. Moreover, for every vertex $i$, one adds
a path $(i|i)$ of length zero, joining $i$ to itself.
A quiver $Q$ now defines an additive category $\catQ$ where all
objects are direct sums of indecomposables and the indecomposable
objects are the elements of $V$, and, for $i,j\in V$,
$\mathop{\rm Mor}_{\catq}(i,j)=$ Set of paths from $i$ to $j$.
\par
Given an abelian category $\catA$, a \it representation of $Q$ in
$\catA$ \rm is a covariant additive functor from $\catQ$ to $\catA$.
Two representations $\catR$
and $\catR^\p$
are called \it equivalent\rm,
if they are isomorphic as functors.  
We denote by $\Rep_{\cata}(Q)$ the abelian category of all representations
of $Q$ in $\catA$.
Notions like \it sub-representations\rm, \it quotient representations\rm,
\it direct sums of representations\rm, etc., are then defined in the usual
way.
\begin{Rem}
A representation $\catR$ of $Q$ in $\catA$ is specified by a collection
$\E_i$, $i\in V$, of objects in $\catA$ and a collection of morphisms
$\phi_a\in \mathop{\rm Mor}_{\cata}(\E_{t(a)},\E_{h(a)})$, $a\in A$.
We then write simply $\catR=(\E_i,i\in V;\phi_a,a\in A)$.
Note that two representations $\catR=(\E_i,i\in V;\phi_a,a\in A)$
and $\catR^\p=(\E^\p_i,i\in V;\phi^\p_a,a\in A)$ are equivalent, if and only if
there are isomorphisms $\psi_i\colon \E_i\lra \E_i^\p$, $i\in V$, with
$\phi_a^\p=\psi_{h(a)}\circ\phi_a\circ \psi_{t(a)}^{-1}$ 
for all arrows $a\in A$.
\end{Rem}
\subsubsection*{Semistability}
Let $(G,\le)$ be a totally ordered abelian group
and $\vartheta\colon {\mathop{\rm Ob}}(\Rep_{\cata}(Q))\lra G$ a map which
factorizes over a group homomorphism $K_{\cata}(Q)\lra G$.
Here, $K_{\cata}(Q)$ is the $K$-group of the abelian category $\Rep_{\cata}(Q)$.
\par
Then, a representation $\catR$ is called \it $\vartheta$-(semi)stable\rm, if and
only if the following two conditions are satisfied
\begin{enumerate}
\item $\vartheta(\catR)=0$
\item $\vartheta(\catR^\p)\q (\le)\q 0$ for every non-trivial proper
      sub-representation $\catR^\p$ of $\catR$.
\end{enumerate}
Recall that "$(\le)$" means that "$<$" is used for defining "stable"
and "$\le$" for defining "semi\-stable".
\subsection*{Representations of quivers in Algebraic Geometry}
Let $Q$ be a quiver and $(X,\O_X(1))$ a polarized 
projective manifold defined over $\C$.
In this paper, we will consider representations of $Q$ in the category
$\Coh(X)$ of coherent sheaves on $X$.
\begin{Conv}
In the following, the word "representation" refers to a representation in
$\Coh(X)$.
\end{Conv}
Let $\catR:=(\E_i, i\in V; \phi_a, a\in A)$ be a representation of $Q$.
The function
\begin{eqnarray*}
\catP\colon V&\lra & \Q[x]\\
             i&\lma & P(\E_i),
\end{eqnarray*}
is called the \it type of $\catR$\rm.
\par
Let $\catP\colon V\lra \Q[x]$ be a map and
$S$ a noetherian scheme. A \it family of representations
of type $\catP$ of $Q$ parametrized by $S$ \rm is
a tuple $({\frak E}_{S,i}, i\in V; \phi_{S,a}, a\in A)$
consisting of $S$-flat families ${\frak E}_{S,i}$
of coherent sheaves with Hilbert polynomials
$\catP(i)$ on $S\times X$, $i\in V$, and 
elements
$\phi_{S,a}
\in \Hom\bigl({\frak E}_{S,t(a)}, {\frak E}_{S, h(a)}\bigr)$,
$a\in A$.
We leave it to the reader to define \it equivalence of families\rm . 
\subsubsection*{Semistability}
\begin{Conv}
From now on, we assume that $Q$ is an oriented tree.
\end{Conv}
We have to find a "good" function $\vartheta$ for which we can prove 
explicit results. For this, the totally order abelian group will be
$(\Q[x],\le)$ where "$\le$" is the lexicographic order of polynomials.
We will now give the definition of $\vartheta$ and explain at the end
of this section why it is the natural
choice.
The definition of $\vartheta$ depends on several parameters, namely,
\begin{itemize}
\item a function $\catP\colon V\lra \Q[x]\setminus\{0\}$,
\item a collection $\ul{\sigma}_Q=(\sigma_i,i\in V)$ 
      of positive rational polynomials $\sigma_i$, $i\in V$,
      of degree at most $\dim X-1$, and
\item a collection $\ul{b}_Q=(b_a,a\in A)$ 
      of positive rational numbers $b_a$, $a\in A$.
\end{itemize}
Having fixed these data, we write $P_i:=\catP(i)$, and let $r_i$ be the
associated rank, $i\in V$. We also set 
$\sigma:=\sigma_1\cdot...\cdot\sigma_n$ and $\check{\sigma}_i
:=\sigma/\sigma_i$, $i\in V$. 
\par
Now define $\vartheta:=\vartheta(\catP,\ul{\sigma}_Q,\ul{b}_Q)$ as the function
which assigns to a representation $\catR=(\E_i,i\in V;\phi_a, a\in A)$
the polynomial
\begin{eqnarray*}
&\sum_{a\in A}b_a&\biggl[\check{\sigma}_{t(a)}\Bigl\{ 
P(\E_{t(a)})-
\rk\E_{t(a)}{P_{t(a)}-\sigma_{t(a)}\over r_{t(a)}}\Bigr\}
\\
&\phantom{\sum_{a\in A}}&\q+\check{\sigma}_{h(a)}\Bigl\{
P(\E_{h(a)})-\rk\E_{h(a)}{P_{h(a)}+\sigma_{h(a)}\over r_{h(a)}}\Bigr\}\biggr].
\end{eqnarray*}
\subsubsection*{Properties of $\vartheta$-semistable representations}
We the above definition of $\vartheta$, we have the concept of 
$\vartheta$-(semi)stability for representations of $Q$ at hand.
Here is a list
of properties of this concept.
In the following $\catR:=(\E_i,i\in V;\phi_a, a\in A)$ is assumed to be
a $\vartheta$-semistable representation, unless otherwise mentioned.
\begin{enumerate}
\item The condition $\vartheta(\catR)=0$ is automatic if the type of $\catR$
      is $\catP$.
\item The sheaves $\E_i$ must all be torsion free. Indeed, set
      $\F_i:=\mathop{\rm Tors}(\E_i)$, $i\in V$, and $\phi_a^\p:=
      \phi_{a|\F_{t(a)}}$, $a\in A$. Then, 
      $\catR^\p:=(\F_i,i\in V;\phi^\p_a, a\in A)$ is a sub-representation 
      of $\catR$ with
      $$
      \vartheta(\catR^\p)\q=\q \sum_{a\in A} b_a \bigl[ \check{\sigma}_{t(a)}
      P(\F_{t(a)})+\check{\sigma}_{h(a)}P(\F_{h(a)})\bigr].
      $$
      This polynomial is strictly positive as soon as one of the sheaves
      $\F_i$ is non-trivial.
\item Suppose the type of $\catR$ is $\catP$. 
      Then all homomorphisms $\phi_a$, $a\in A$,
      are non-zero. If, say, $\phi_{a_0}$ were zero, we could remove the arrow
      $a_0$ from $Q$ in order to obtain two disjoint subtrees $Q_{t(a_0)}$ and
      $Q_{h(a_0)}$. We define $\F_i$ as $\E_i$ if $i\in V_{t(a_0)}$ and as
      $0$ otherwise, and $\phi_a^\p:=\phi_{a|\F_{t(a)}}$ for all $a\in A$.
      Thus, $\catR^\p=(\F_i,i\in V; \phi_a^\p, a\in A)$ is a 
      sub-representation
      of $\catR$ with
      $$
      \vartheta(\catR^\p)\q=\q b_{a_0}\cdot\sigma\cdot r_{t(a_0)}\q>\q0,
      $$
      contradicting the assumption that $\catR$ be $\vartheta$-semistable. 
\item The representation $\catR$ possesses a 
      \it Jordan-H\"older filtration\rm
      $$
      0=:\catR^{(m+1)}\subset \catR^{(m)}\subset \cdots \subset \catR^{(1)}
      \subset \catR^{(0)}:=\catR
      $$
      where $\catR^{(l)}$ is a sub-representation of $\catR^{(l-1)}$
      which is maximal w.r.t.\ inclusion among those sub-representations
      $\catR^\p$ with $\vartheta(\catR^\p)=0$, $l=1,...,m+1$.
      The successive quotients $\catR^{(l)}/\catR^{(l+1)}$, $l=0,...,m$, are thus
      $\vartheta$-stable representations and the \it associated graded
      object\rm
      $$
      \gr(\catR)\q:=\q \bigoplus_{l=0}^m \catR^{(l)}/\catR^{(l+1)}
      $$
      which is well-defined up to equivalence is again a $\vartheta$-semistable
      representation of the same type as $\catR$.
      \par
      As usual, two $\vartheta$-semistable representations $\catR$ and $\catR^\p$
      are called \it S-equivalent\rm, if their associated graded objects
      are equivalent, and $\catR$ is called \it $\vartheta$-polystable\rm, 
      if it is
      equivalent to $\gr(\catR)$.
\item One can "join" $\vartheta$-semistable representations:
      Let $Q_1$ and $Q_2$ be subquivers of $Q$, such that $Q=Q_1\cup Q_2$ and 
      $A_1\cap A_2=\varnothing$, 
      and $\catR=(\E_i,i\in V;\phi_a,a\in A)$ a 
      --- not necessarily $\vartheta$-semistable --- representation
      of $Q$. Suppose the representation $\catR_j:=(\E_i,i\in V_j, \phi_a,
      a\in A_j)$ is 
      $\vartheta(\catP_j,\ul{\sigma}_{Q_j},\ul{b}_{Q_j})$-(semi)stable
      for $j=1,2$. Then, $\catR$ is $\vartheta$-(semi)stable. Here, the data
      $(\catP_j,\ul{\sigma}_{Q_j},\ul{b}_{Q_j})$ are obtained from
      $(\catP,\ul{\sigma}_Q,\ul{b}_Q)$ by restriction to $Q_j$, $j=1,2$.
\end{enumerate}
\begin{Ex}[Holomorphic triples]
A \it holomorphic triple \rm is a triple $(\E_1,\E_2,\allowbreak\phi)$
consisting of two coherent sheaves $\E_1$ and $\E_2$
and a homomorphism $\phi\colon \E_1\lra\E_2$.
In other words, a holomorphic triple is a representation
of the quiver $\bullet\lra\bullet$.
Specializing our general definitions to holomorphic triples,
we say that --- for given positive polynomials $\sigma_1$ and
$\sigma_2\in\Q[x]$ of degree at most $\dim X-1$ ---
a holomorphic triple $(\E_1,\E_2\,\phi)$ is  
$\vartheta$-(semi)stable,
if for any two
subsheaves $\F_1$ and $\F_2$ of $\E_1$ and $\E_2$, respectively,
such that $0\neq \F_1\oplus\F_2\neq \E_1\oplus\E_2$ and
$\phi(\F_1)\subset\F_2$
$$
\sigma_2\left(P(\F_1)-\rk\F_1\left({P(\E_1)\over\rk\E_1}-
{\sigma_1\over\rk\E_1}\right)\right)
+ 
\sigma_1\left(P(\F_2)-\rk\F_2\left({P(\E_2)\over\rk\E_2}+
{\sigma_2\over\rk\E_2}\right)\right)
$$
is a (non-positive) negative polynomial in the lexicographic order
of polynomials. Here, $\vartheta$ is associated with $i\lma P(\E_i)$,
$\sigma_1$, $\sigma_2$, and $b_a=1$.
\par
If $X$ is a curve and $\sigma_1=\sigma_2=:\sigma\in\Q_+$,
set $\tau:=\mu(\E_2)+\sigma/\rk \E_2$.
Then the above definition yields the definition of
$\tau$-(semi)stability of Bradlow and Garcia-Prada \cite{BG},
\cite{GP} for holomorphic triples.
Thus, we see that our concept of semistability is more general
because it involves two parameters instead of one.
This might be useful for comparing the moduli spaces
of holomorphic triples with other moduli spaces. 
\end{Ex}
\begin{Rem}  
Let $\catR=(\E_i,i\in V;\phi_a, a\in A)$
be a representation of type $\catP$. By property 5.,
the $\vartheta$-(semi)stability condition is satisfied,
if, for every arrow $a\in A$, $(\E_{t(a)},\E_{h(a)},
\allowbreak \phi_a)$ 
is a 
$\vartheta_a$-(semi)stable holomorphic triple, $\vartheta_a$ being obtained
by restricting the data $(\catP, \ul{\sigma}_Q, \ul{b}_Q)$ to the
subquiver $t(a)\lra h(a)$, $a\in A$.
Therefore, if all $\sigma_i$'s are equal to some $\sigma$, 
existence theorems for $\vartheta$-(semi)stable
representations on curves can be extracted
from the work of Bradlow and Garcia-Prada \cite{BG}. 
\end{Rem}
\begin{Ex}
\label{king1}
Assume that our base manifold $X$ is just a point. Then, $\Coh(X)$ is the
category of finite dimensional complex vector spaces, i.e.,
a representation is of the form $(E_i,i\in V; f_a, a\in A)$ where
$E_i$ is a finite dimensional $\C$-vector space, $i\in V$, and
$f_a\colon E_{t(a)}\lra E_{h(a)}$ is a linear map, $a\in A$.
\par
In this case, the datum $\ul{\sigma}_Q$ is obsolete (strictly speaking,
not defined). This means, we just fix $\catP\colon V\lra \Z_{>0}$
and $\ul{b}_Q=(b_a, a\in A)$ and define
$$
\vartheta(E_i,i\in V; f_a, a\in A)\q:=\q
\sum_{a\in A} b_a\left(\frac{\dim E_{t(a)}}{\catP(t(a))}-
                       \frac{\dim E_{h(a)}}{\catP(h(a))}\right).
$$
The corresponding concept of $\vartheta$-(semi)stability 
agrees with King's notion of $\chi$-(semi)stab\-ility 
for representations $\catR$ with
dimension vector (=type) $\catP$ associated to the character
$$
\begin{array}{rccc}
\chi\colon & \GL(\catP(1))\times\cdots\times \GL(\catP(n)) &\lra& \C^*\\
           & (m_1,...,m_n)&\lma& \det(m_1)^{s_1}\cdot...\cdot \det(m_n)^{s_n},
\end{array}
$$
with
$$
s_i\q:=\q \sum_{a:t(a)=i} \frac{b_a}{\catP(i)}- 
          \sum_{a:h(a)=i} \frac{b_a}{\catP(i)},\qquad i=1,...,n.
$$
This shows that our definition is a natural extension of King's (specialized
to oriented trees) to higher
dimensions.
\end{Ex}      
\subsection*{The main result}
%
Fix $\catP$, $\ul{\sigma}_Q$, and $\ul{b}_Q$ 
as before, and set $\vartheta:=\vartheta(\catP,\ul{\sigma}_Q,\ul{b}_Q)$.
Define $\catM(Q)_{\catp}^{\vartheta-(s)s}$                 
as the functor which assigns to a noetherian scheme $S$ the set
of equivalence classes of families of 
$\vartheta$-(semi)stable  
representations of type $\catP$ of $Q$ which are parametrized by $S$.
\begin{Thm}
\label{MTH}
{\rm i)} There exist a projective scheme
${\cal M}:={\cal M}(Q)_{\catp}^{\vartheta-ss}$
and a natural transformation 
$\catT\colon \catM(Q)_{\catp}^{\vartheta-ss}
\lra h_{\cal M}$, such that for any other scheme ${\cal M}^\p$
and any natural transformation $\catT^\p$, there exists 
a unique morphism $\rho\colon {\cal M}\lra {\cal M^\p}$
with $\catT^\p=h(\rho)\circ \catT$.
\par
{\rm ii)} The map $\catT(\pt)$ induces a bijection between
the set of S-equivalence classes of $\vartheta$-semistable
representations of $Q$ of type $\catP$ and the set
of closed points of ${\cal M}$.
\par
{\rm iii)} The space ${\cal M}$
contains an open subscheme ${\cal M}(Q)
_{\catp}^{\vartheta-s}$
which becomes through $\catT$ a coarse moduli scheme for 
the functor
$\catM(Q)_{\catp}^{\vartheta-s}$.
\end{Thm}
\begin{Ex}
\label{king2}
As an illustration how such a moduli space is constructed and as a tool
for later sections, we review the case $X=\{\pt\}$, i.e., King's 
construction
in the case $Q$ is an oriented tree.
\par
Fix $\catP\colon V\lra \Z_{>0}$, $\ul{b}_Q$, and write
$\vartheta:=\vartheta(\catP,\ul{b}_Q)$, and let $\chi$ be the corresponding
character (see \ref{king1}) of
$$
\GL(\catP)\q:=\q \GL(\catP(1))\times\cdots\times \GL(\catP(n)).
$$
Every representation of $Q$ with dimension vector $\catP$ is equivalent
to one in the space
$$
{\frak H}(\catP)
\q :=\q
\bigoplus_{a\in A} \Hom\bigl(\C^{\catp(t(a))}, \C^{\catp(h(a))}\bigr).
$$
On this space, there is a natural left action of
$\GL(\catP)$, and
the set of equivalence classes of representations with dimension
vector $\catP$ corresponds to the set of $\GL(\catP)$-orbits in
${\frak H}(\catP)$. However, this set does not carry a natural structure
of an algebraic variety or even of a topological space.
\par
In order to obtain an algebraic variety, we must use the GIT machinery,
i.e., we must choose a linearization of the given action in 
$\O_{{\frak H}(\catp)}$. This is given by the character $\chi$.
King's moduli space is the GIT quotient 
${\frak H}(\catP)\catqot_\chi\GL(\catP)$. Note that one must show that the
points which are (semi)stable w.r.t.\ the given linearization
are exactly the ones corresponding to $\chi$-(semi)stable representations. 
\par
We can describe the quotient in another way. 
Note that we know (Property 3.) that every homomorphism $f_a$
occuring in a $\vartheta$-semistable representation $(E_i,i\in V;f_a,a\in A)$
must be non-zero. Moreover, one has
\begin{Lem}
\label{Mult}
Let $Q$ be an oriented  tree and
$(f_a, a\in A)$ be a point in ${\frak H}(\catP)$.
Then, for a given $a_0\in A$ and $z\in\C^*$, there exists an element
$g_{a_0,z}=(z_1\id,...,\allowbreak z_n\id)\in\GL(\catP)$ such that
$(f_a^\p, a\in A):= g_{a_0,z}\cdot(f_a, a\in A)$
looks as follows: $f_a^\p=f_a$ for $a\neq a_0$,
and $f_{a_0}^\p=z\cdot f_{a_0}$.
\end{Lem}
{\it Proof}.
Removing the arrow $a_0$ from the quiver $Q$ yields two connected
subquivers $Q_{t(a_0)}$ and $Q_{h(a_0)}$.
We set $z_i:=z$ for $i\in V_{h(a_0)}$ and $=1$ for $i\in V_{t(a_0)}$.
Then, $g_{a_0,z}:=(z_1\id,...,z_n\id)$ does the trick.
\par
\medskip
Thus, we can start with the space
$$
{\frak P}(\catP)\q :=\q \prod_{a\in A}\P\left(\Hom\bigl(\C^{\catp(t(a))}, 
\C^{\catp(h(a))}\bigr)^\vee\right)
$$
with fixed ample line bundle
$\O(b_a,a\in A)$.
The moduli space is now the projective variety
${\frak P}(\catP)\catqot_{\ul{b}_Q} \SL(\catP)$.
Here, it will follow from Theorem~\ref{GITII} that the 
$\O(b_a,a\in A)$-(semi)stable
points are exactly those which correspond to $\vartheta$-(semi)stable
representations.
It is easy to see that this is the same variety as 
${\frak H}(\catP)\catqot_\chi\GL(\catP)$.
\par
Now, let $X$ be an arbitrary projective manifold.
The GIT construction in this case follows the same pattern:
\begin{itemize}
\item Step 1: Find a variety ${\frak T}$ analogous to 
              ${\frak P}(\catP)$ which
              parametrizes representations of type $\catP$ and contains
              every $\vartheta$-semistable representation at least
              once.
\item Step 2: Show that there is an action of a reductive algebraic
              group $G$, such that two points in ${\frak T}$
              lie in the same orbit if and only if they correspond to
              equivalent representations.
\item Step 3: Find a linearization, such that a point in ${\frak T}$
              is (semi)stable w.r.t.\ that linearization if and only if
              it corresponds to a $\vartheta$-(semi)stable representation.
\end{itemize}
Having treated all these steps successfully, one gets the moduli space
as the GIT quotient ${\frak T}\catqot G$.
In our construction, the assumption that $Q$ be an oriented tree is essential.
First, it allows us as in the example to choose $G$ as a product of
special linear groups. This makes the computations for the Hilbert-Mumford
criterion already simpler. Second, for the action of $\SL(\catP)$
on ${\frak P}(\catP)$, one has Theorem~\ref{GITII} which simplifies the
computations even further. But this theorem is also crucial for
proving that one can in fact adjust all parameters appearing in such a way
that Step 1 - 3 really go through.
\end{Ex}
\subsection*{Concluding Remarks}
\subsubsection*{How did we find $\vartheta$?}
The most general moduli problem one would like to treat is the following:
Let $G$ be a reductive algebraic group and $Y$ a projective manifold
on which $G$ acts. If $P$ is a principal $G$-bundle on $X$, we obtain
an induced fibre space $Y(P):=P\times^G Y$. One would now like to
classify pairs $(P,\sigma)$ where $P$ is a principal $G$-bundle
and $\sigma\colon X\lra Y(P)$ is a section.
For this, one has to define a general semistability concept and
establish the existence of moduli spaces.
In gauge theory, such a programme has been succesfully treated in
\cite{Ba} and \cite{MR}.
In the algebraic context, the author has defined this semistability
concept and constructed the moduli spaces in the case when $X$ is a curve,
$G=\GL(r)$, and
the action of the center  $\C^*\cdot E_n\subset G$ 
on $Y$ is trivial. The semistability
concept is a version of the Hilbert-Mumford criterion and depends only
on the choice of a linearization of the $G$-action on $Y$ \cite{Sch3}.
This concept is completely natural and reproduces all known examples.
\par
In the present paper, we have $G=\GL(\catP)$ and $Y={\frak P}(\catP)$. Note
that the induced ${\C^*}^n$-action on $Y$ is trivial.
One can adapt the techniques of \cite{Sch3} to the present situation.
We illustrate this by an example:
\par
Let us look at the case when $X$ is a curve
and $Q=1\lra 2\lra 3$.
The natural parameter space for a bounded family of triples
of vector bundles $(E_1,E_2,E_3)$ where $\deg E_i=:d_i$ and
$\rk E_i=:r_i$ are fixed, $i=1,2,3$, is a product of quot-schemes
${\frak Q}:={\frak Q}_1\times {\frak Q}_2\times {\frak Q}_3$.
If we polarize ${\frak Q}_i$ by the line bundle $\O_i(1)$ coming from 
Gieseker's covariant map, we can polarize
${\frak Q}$ by $\O(\tau_1,\tau_2,\tau_3)$
where the $\tau_i$ are positive rational numbers.
We choose in addition a polarization $\O(b_1,b_2)$ on
$\P(\Hom(\C^{r_1},\C^{r_2})^\vee)\times \P(\Hom(\C^{r_2},\C^{r_3})^\vee)$.
\par
Let $\catR :=(E_1,E_2,E_3;\phi_1,\phi_2)$ be a representation of $Q$.
As explained in \cite{Sch3}, the objects for testing semistability
are triples 
$$
\ul{T}=\bigl((\ul{E}^{1,\bullet},\ul{\alpha}_1),
(\ul{E}^{2,\bullet},\ul{\alpha}_2),(\ul{E}^{3,\bullet},\ul{\alpha}_3)\bigr), 
$$
where
$$
(\ul{E}^{i,\bullet},\ul{\alpha}_i)=\bigl(0\subset E^{i,1}\subset \cdots \subset E^{i, s_i}\subset E,
(\alpha_{i,1},...,\alpha_{i,s_i})\bigr)
$$ 
is a weighted filtration for $E_i$, $i=1,2,3$.
Recall that
$$
M\bigl(\ul{E}^{i,\bullet},\ul{\alpha}_i\bigr)
\q:=\q
\sum_{j=1}^{s_i}\alpha_{i,j}\bigl(\deg E_i\rk E^{i,j}-\deg E^{i,j} r_i\bigr),\q i=1,2,3.
$$ 
One defines $\mu_{b_1,b_2}(\ul{T};\phi_1,\phi_2)$ similarly as in \cite{Sch3}.
\par
With these conventions, $\catR$ is called \it 
$(\tau_1,\tau_2,\tau_3;b_1,b_2)$-(semi)stable\rm, if for every triple of weighted
filtrations $\ul{T}$ as above one finds
$$
\sum_{i=1}^3 \tau_i M\bigl(\ul{E}^{i,\bullet},\ul{\alpha}_i\bigr)
+
\mu_{b_1,b_2}\bigl(\ul{T};\phi_1,\phi_2\bigl)\q(\ge)\q 0.
$$
\par
Now, the decomposition results of Section~\ref{GITI} permit us to restrict
to triples $\ul{T}$ where  
$$
(\ul{E}^{i,\bullet},\ul{\alpha}_i)\q=\q 
\bigl(0\subset F_i\subset E, (1/r_i)\bigr)
$$
for some subbundles $F_i$ of $E_i$, $i=1,2,3$, with $\phi_1(F_1)\subset F_2$
and $\phi_2(F_2)\subset F_3$.
In this case, one has
$$
\mu_{b_1,b_2}\bigl(\ul{T};\phi_1,\phi_2\bigr)\q=\q
b_1\left(\frac{\rk F_2}{r_2}-\frac{\rk F_1}{r_1}\right)+b_2\left(
\frac{\rk F_3}{r_3}-\frac{\rk F_2}{r_2}\right),
$$
so that the condition becomes
\begin{eqnarray*}
\tau_1\bigl(\rk F_1 \mu(E_1)-\deg F_1\bigr)
+\tau_2 \bigl(\rk F_2 \mu(E_2)-\deg F_2\bigr)
+\tau_3 \bigl(\rk F_3 \mu(E_3)-\deg F_3\bigr) 
\\
+b_1\left(\frac{\rk F_2}{r_2}-\frac{\rk F_1}{r_1}\right)+b_2\left(
\frac{\rk F_3}{r_3}-\frac{\rk F_2}{r_2}\right)\q(\ge)\q0.&&
\end{eqnarray*}
Now, write $\tau_1=b_1/\sigma_1$, $\tau_2=(b_1+b_2)/\sigma_2$,
and $\tau_3=b_2/\sigma_3$ for some positive rational numbers 
$\sigma_i$. We thus see that $\catR$ will be 
$(\tau_1,\tau_2,\tau_3;b_1,b_2)$-(semi)stable, if and only if for every
sub-representation $(F_1,F_2,F_3;\phi_1^\p,\phi_2^\p)$,
such that $F_1\oplus F_2\oplus F_3$ is a non-trivial proper subbundle
of $E_1\oplus E_2\oplus E_3$, one has
\begin{eqnarray*}
b_1\left(\frac{1}{\sigma_1}\Bigl(\rk F_1\frac{d_1-\sigma_1}{r_1}-\deg F_1\Bigr)
+ \frac{1}{\sigma_2}\Bigl(\rk F_2\frac{d_1+\sigma_2}{r_2}-\deg F_2\Bigr)\right)
&+&
\\
b_2\left(\frac{1}{\sigma_2}\Bigl(\rk F_2\frac{d_2-\sigma_2}{r_2}-\deg F_2\Bigr)
+ \frac{1}{\sigma_3}\Bigl(\rk F_3\frac{d_3+\sigma_3}{r_3}-\deg F_3\Bigr)\right)
&(\ge)& 0.
\end{eqnarray*}
Multiply this by $-\sigma_1\sigma_2\sigma_3$ to recover our original definition.
\subsubsection*{Other quivers}
One would expect to be able to treat other quivers as well, at least when
they don't contain oriented cycles.
In this case, the concept of semistability should depend on $\catP$ and
$\ul{\sigma}_Q$ as before and a character $\chi$ of the group
$\GL(\catP(1))\times \cdots\times \GL(\catP(n))$, but I expect it to look
much more difficult.
\par
For quivers with oriented cycles, there arise other problems.
Look for example at the theory
of Higgs bundles which can
be viewed as the moduli problem associated with the quiver
consisting of one vertex and an arrow, joining the vertex to itself.
Applying the above definition of a representation would lead
to the consideration of sheaves $\E$ together with an
endomorphism $\phi\colon\E\lra\E$. It turns out that the 
(semi)stability concept forces $\E$ to be a (semi)stable sheaf
(see \cite{Sch}, Thm.~3.1, with $L=\O_X$).
As, moreover, a stable sheaf has no endomorphisms besides multiples
of the identity, this theory brings nothing new.
The way out is to consider twisted endomorphisms
$\phi\colon \E\lra\E\otimes L$, $L$ a suitable line bundle, e.g., 
$L=K_X$, if $X$ is a curve of genus $g\ge 2$ \cite{Hi}.
\par
Finally, the resulting moduli spaces are only quasi-projective.
In order to compactify them, one must add further data \cite{Sch},
\cite{Sch3}. If one does this, one finds also semistability concepts
which cannot be formulated as conditions on sub-representations only
\cite{Sch3}.
\section{More notation concerning quivers}
We introduce now some terminology for quivers which is adapted to
the subsequent proofs.
\par
Let $Q=(V,A,t,h)$ be a quiver.
A \it subquiver \rm $Q^\p\subset Q$ consists of
a subset $V^\p\subset V$ and a subset $A^\p\subset A$, such that
$h(A^\p)\cup t(A^\p)\subset V^\p$, and
is called a \it full subquiver\rm, if 
any arrow $a\in A$ with $h(a)\in V^\p$ and $t(a)\in V^\p$ lies
in $A^\p$. Obviously, we can associate
to any subset $V^\p\subset V$ a full subquiver $Q(V^\p)$ of $Q$.
For any $i\in V$, 
\it the star of $i$ \rm is defined as the subquiver
$\St_Q(i)$ of $Q$ whose set of arrows is $A(i):=\{\, a\in A\,|\, h(a)=i
\ \vee\ t(a)=i\,\}$ and whose set of vertices is $V(i):=h(A(i))\cup t(A(i))$.
A vertex $i$ of the subquiver 
$Q^\p$ is called an \it end (of $Q^\p$ in $Q$)\rm,
if $\St_Q(i)$ is not contained in $Q^\p$. The set of all ends will be
denoted by $\mathop{\rm END}_Q(Q^\p)$. For each vertex 
$i\in \mathop{\rm END}_Q(Q^\p)$, the set of \it ingoing arrows \rm
$\mathop{\rm In}_Q(i)$ is defined as the set of all arrows
of $\St_Q(i)$ not lying in $Q^\p$ whose head is $i$.
Similarly, we define $\mathop{\rm Out}_Q(i)$, the set of
\it outgoing arrows\rm .
\par
The following lemma will enable us to prove many of the needed
technical details by induction.
\begin{Lem}
\label{Cutoff}
Let $Q$ be an oriented  tree, then there exists a vertex $i$
whose star is either $\bigl(\{i,i^\p\}, (i,i^\p)\bigr)$
or $\bigl(\{i,i^\p\}, (i^\p,i)\bigr)$ for some vertex $i^\p\in V$.
\end{Lem}
So, after relabelling the vertices, we can assume that $i=n$ and $i^\p=n-1$
and define a new quiver $Q^\p$ with $V^\p:=\{\, 1,...,n-1\,\}$
and $A^\p=A\setminus\{(n-1,n)\}$ or
$A^\p=A\setminus\{(n,n-1)\}$. This is again
an oriented  tree.
\section{Decomposition of one parameter subgroups}
\label{GITI}
In this section, we prove the main auxiliary result which simplifies
the Hilbert-Mumford criterion for the actions we consider.
At a first reading, the reader might follow this Section till
Theorem~\ref{GITII} and then proceed directly to the proof of the main
result.
\par
Let $Q$ be an oriented  tree with $V=\{\, 1,...,n\,\}$.
Let $V_1,...,V_n$, and $W_1,...,W_n$ 
be finite dimensional $\C$-vector spaces
and suppose we are given representations
$\tau_i\colon \SL(V_i)\lra\GL(W_i)$, $i=1,...,n$.
Set $\ul{\tau}:=(\tau_1,...,\tau_n)$, define $\catP\colon V\lra \Z_{>0}$,
$i\lma p_i:=\dim V_i$, and ${\frak P}(\catP):=\prod_{a\in A}
\P(\Hom(V_{t(a)}, V_{h(a)})^\vee)$.
These data define an action of $\SL(\catP):=\prod_{i\in V}\SL(V_i)$ on
$$
\P_{\ul{\tau},\catp}\q:=\q \P(W_1^\vee)\times\cdots\times\P(W_n^\vee)
\times {\frak P}(\catP).
$$
Fix a (fractional) 
polarization $\O(l_1,...,l_n; b_a, a\in A)$ on the space
$\P_{\ul{\tau},\catp}$ where the $l_i$ and $b_a$
are positive rational numbers. 
It will be our task to describe the (semi)stable points in
$\P_{\ul{\tau},\catp}$ w.r.t.\ the given linearization.
\subsection*{Further assumptions and notations}
A one parameter subgroup of $\SL(\catP)$ will be written as 
$\la=(\la_1,....,\la_n)$ where $\la_i$ is a one parameter subgroup
of $\SL(V_i)$, $i=1,...,n$.
Let $\ul{w}:=([w_1],...,\allowbreak [w_n]; [f_a], a\in A)$ be a point
in $\P_{\ul{\tau},\catp}$ 
and $\la$ be a one parameter subgroup of $\SL(\catP)$.
Then, $\mu(\ul{w},\la)$ is defined as minus the weight of
the induced $\C^*$-action on the fibre of $\O(l_1,...,l_n,b_a, a\in A)$
over $\ul{w}_\infty:=\lim_{z\ra\infty} \la(z)\cdot \ul{w}$.
Recall that the Hilbert-Mumford critrion states that $\ul{w}$ is (semi)stable
w.r.t.\ the linearization in $\O(l_1,...,\allowbreak
l_n; b_a, a\in A)$ if and only if
$\mu(\ul{w},\la)\ (\ge)\ 0$ for all non trivial one parameter subgroups
$\la$ of $\SL(\catP)$.
\par
Now, let $\ul{w}$ and $\la$ be as before.
For $i=1,...,n$, write $\mu([w_i], \la_i)$ for the weight
of the $\C^*$-action induced by the action of $\la_i$ on the fibre
of $\O_{\P(W_i^\vee)}(-1)$ over the point $\lim_{z\ra\infty}\la_i(z)[w_i]$,
and, for an arrow $a\in A$, we let $\mu([f_a], (\la_{t(a)},\la_{h(a)}))$
be the weight of the resulting $\C^*$-action on the fibre 
of $\O_{\P(\Hom(V_{t(a)},V_{h(a)})^\vee)}(-1)$ over 
$\lim_{z\ra \infty} (\la_{t(a)}(z),\la_{h(a)}(z))\cdot [f_a]$.
With these conventions
$$
\mu(\ul{w},\la)\q=\q l_1\mu([w_1],\la_1)+\cdots+l_n\mu([w_n],\la_n)
+\sum_{a\in A}b_a\mu\bigl([f_a], (\la_{t(a)},\la_{h(a)})\bigr).
$$
Next, we remind you that a one parameter subgroup $\la_i$
of $\SL(V_i)$ is defined by giving a basis $v_1^i,...,v_{p_i}^i$
for $V_i$ and integer weights $\gamma_1^i\le\cdots\le\gamma_{p_i}^i$
with $\sum_j \gamma_j^i=0$.
We will also use \it formal one parameter subgroups \rm
which are defined as before, only that this time
the $\gamma^i_j$ are allowed to be rational numbers.
It is clear how to define $\mu(\ul{w}, \la)$
for a formal one parameter subgroup.
Let $\la=(\la_1,...,\la_n)$ be a (formal) one parameter subgroup
where $\la_i$ is given (w.r.t.\ to some basis of $V_i$) by the
weight vector $\ul{\gamma}^i$. Then, we write $\ul{\gamma}=
(\ul{\gamma}^1,...,\ul{\gamma}^n)$.
For $j=0,...,p_i$, we can look at the one parameter subgroup
$\la_i^{(j)}$ which is defined by the weight vector
$\gamma^{i,(j)}:=(\, j-p_i,...,j-p_i, j,...,j\,)$, $j-p_i$ 
occuring $j$ times. Note that both $\la_i^{(0)}$ and $\la_i^{(p_i)}$ are
the trivial one parameter subgroup.
The $\la_i^{(j)}$ with $1\le j<p_i$ are particularly important due to the fact that any
weight vector $\ul{\gamma}^i=(\, \gamma_1^i,...,\gamma_{p_i}^i\,)$ as before
can be written as 
\begin{equation}
\label{Decomp0}
\ul{\gamma}^i\q =\q \sum_{j=1}^{p_i-1}
{(\gamma^i_{j+1}-\gamma^i_j)\over p_i}\gamma^{i,(j)}.
\end{equation}
\begin{Ass}
\label{Add}
We require the following additivity property for the action
of $\SL(V_i)$ on $\P(W_i^\vee)$, $i=1,...,n$: 
For every point $[w_i]\in \P(W_i^\vee)$, every basis
$v^i_1,...,v^i_p$ of $W_i$, and every two one parameter subgroups $\la$
and $\la^\p$ of $\SL(V_i)$ which are given with respect to that basis
by weight vectors $(\gamma_1,...,\gamma_p)$ and 
$(\gamma_1^\p,...,\gamma_p^\p)$
with $\gamma_1\le \cdots\le \gamma_p$ and 
$\gamma_1^\p\le\cdots\le\gamma_p^\p$,
we have
$$
\mu([w_i],\la\cdot\la^\p)\q=\q 
\mu([w_i],\la)+\mu([w_i],\la^\p).
$$
\end{Ass}
\begin{Ex}
In general, one has in the above situation only ``$\le$'', e.g.,
for the action of $\SL(V)$ on $\End(V)$. To see this, consider
for example $V=\C^3$ and the action of $\SL_3(\C)$ on $M_3(\C)$,
the vector space of complex $(3\times 3)$-matrices, by conjugation.
Let $f$ be given by the matrix
$$
\left(
\begin{array}{ccc}
0 & 1& 0\\
0&0&1\\
0&0&0
\end{array}
\right)
$$
and $\la$ and $\la^\p$ with respect to the standard
basis by $(-2,1,1)$ and $(-1,-1,2)$, respectively.
Then, $\mu([f], \la)=0=\mu([f],\la^\p)$, but 
$\mu([f], \la\cdot\la^\p)=-3$.
A similar phenomenon is responsible for the technicalities
which we will encounter below.
\end{Ex}
If we are given bases $v^i_1,...,v_{p_i}^i$ for $V_i$, we set
$V^{(j)}_i:=\langle\, v_1^i,...,v_j^i\,\rangle$, $j=0,...,p_i$, $i=1,...,n$.
Let $\ul{j}:=(j_1,...,j_n)$ be a tuple of elements
with $j_i\in \{\, 0,...,p_i\,\}$, $i=1,...,n$.
To such a tuple we associate the weight vector
$$
\ul{\gamma}^{\ul{j}}\q:=\q\Bigl(\, {1\over p_1}\gamma^{1,(j_1)},...,
{1\over p_n}\gamma^{n,(j_n)}\,\Bigr),
$$
and denote the corresponding formal one parameter 
subgroup of $\SL(\catP)$ by
$\la^{\ul{j}}$.
\begin{Thm}
\label{GITII}
In the above setting, assume that there are positive rational numbers
$\alpha_{i,a}$, $i\in V$, $a\in A$, with
$$
l_i\q=\q \sum_{a\in A: t(a)=i\vee h(a)=i} b_a \alpha_{i,a},
$$
then a point 
$([w_1],...,[w_n],[f_a],a\in A)$ 
is (semi)stable w.r.t.\ the
linearization of the action of $\SL(\catP)$ on $\P_{\tau,\catp}$
in $\O(l_1,...,l_n;b_a,a\in A)$, 
if and only if
for all possible choices of bases $v^i_1,...,v_{p_i}^i$ for $V_i$, $i\in V$,
and indices
$j_i\in \{\, 0,...,p_i\,\}$ with
$$
f_a\bigl(V_{t(a)}^{(j_{t(a)})}\bigr)\q\subset\q
V_{h(a)}^{(j_{h(a)})},\qquad \hbox{for all } a\in A, 
$$ 
one has
\begin{eqnarray*}
0&(\le)&\sum_{a\in A} b_a\Bigl[\alpha_{t(a),a}{1\over p_{t(a)}}
\mu\bigl([w_{t(a)}], \la_{t(a)}^{(j_{t(a)})}\bigr)-{j_{t(a)}\over p_{t(a)}}
\\
&&
\phantom{\sum_{a\in A} b_a\Bigl[
} + \alpha_{h(a),a}{1\over p_{h(a)}}
\mu\bigl([w_{h(a)}], \la_{h(a)}^{(j_{h(a)})}\bigr)+{j_{h(a)}\over p_{h(a)}}
\Bigr].
\end{eqnarray*}
\end{Thm} 
The rest of this Section concerns the proof of Theorem~\ref{GITII}.
Let $([f_a],a\in A)$ be an element in ${\frak P}(\catP)$.
We call the weight vector $\ul{\gamma}^{\ul{j}}$ \it basic
(w.r.t.\ $([f_a],a\in A)$)\rm, if (1) the subquiver
$Q_{\ul{j}}:=Q(\{\, i\in V\, |\, 0<j_i<p_i\,\})$ is connected,
(2) for any arrow $a\in A_{\ul{j}}$ we have
$f_a\bigl(V^{(j_{t(a)})}_{t(a)})\subset V_{h(a)}^{(j_{h(a)})}$,
and (3) neither $V_{t(a)}^{(j_{t(a)})}\subset\ker f_a$
nor $V_{h(a)}^{(j_{h(a)})}\supset\Im f_a$, $a\in A$.
The strategy is now to decompose the weight vector of
any given one parameter subgroup in a suitable way into basic ones,
so that a point will be (semi)stable if and only if
the Hilbert-Mumford criterion is satisfied for basic formal
one parameter subgroups. 
\subsection*{The central decomposition theorem}
\subsubsection*{The case $\#V=2$}
Let $[f]\in{\frak P}(\catP)$ be the class
of a homomorphism $f\colon V_1\lra V_2$.
Given bases $v_1^j,...,v_{p_j}^j$ for $V_j$, $j=1,2$,
we write $f=\sum_{i,j} f_{i,j} {v_i^1}^\vee
\otimes v_j^2$.
\begin{Thm}
\label{Decomp1}
Let $[f]$ be as before fixed.
Then, for any given one parameter subgroup $(\la_1,\la_2)$ of
$\SL(V_1)\times \SL(V_2)$ which is specified by the bases
$v^1_1,...,v_{p_j}^j$ of $V_j$, $j=1,2$, and the weight vector
$(\ul{\gamma}^1,\ul{\gamma}^2)$, there exist indices $i(1)_*$ and 
$i(2)_*$
with $f_{i(1)_*i(2)_*}\neq 0$ such that $\mu([f],(\la_1,\la_2))=
\mu([{v^1_{i(1)_*}}^\vee
\otimes v^2_{i(2)_*}],(\la_1,\la_2))$ and a decomposition
$$
(\ul{\gamma}^1,\ul{\gamma}^2)\q =\q
\sum_{{i(1)=0,...,p_1-1; \atop
       i(2)=0,...,p_2-1\phantom{;}}} \eta_{i(1)i(2)}\Bigl({1\over p_1}
\gamma^{1,(i(1))}, {1\over p_2} \gamma^{2,(i(2))}\Bigr),
$$
such that all the coefficients are non-negative rational numbers,
and whenever the coefficient
$\eta_{i(1)i(2)}$ is not zero, the weight vector
$(\gamma^{1, (i(1))}, \gamma^{2, (i(2))})$ is a basic weight vector
and $\mu([f], \la^{(i(1), i(2))})=\mu([{v^1_{i(1)_*}}^\vee
\otimes v^2_{i(2)_*}],  \la^{(i(1), i(2))})$.
\end{Thm}
One infers
\begin{eqnarray*}
&&\mu\bigl(([w_1], [w_2], [f]), (\la_1,\la_2)\bigr)\\
 &=&
\mu\bigl(([w_1], [w_2],[{v^1_{i(1)_*}}^\vee
\otimes v^2_{i(2)_*}]),(\la_1,\la_2)\bigr)
\\
&=&
\sum \eta_{i(1)i(2)} \mu\bigl(([w_1], [w_2], [{v^1_{i(1)_*}}^\vee
\otimes v^2_{i(2)_*}]), \la^{(i(1),i(2))}\bigr)
\\
&=&
\sum \eta_{i(1)i(2)} \mu\bigl(([w_1], [w_2], [f]), \la^{(i(1),i(2))}\bigr),
\end{eqnarray*}
whence we have achieved our goal in this situation.
\par
\medskip
\it Proof of Theorem~\ref{Decomp1}\rm.
To reduce indices, we slightly change the notation: We set $V:=V_1$ and
$W:=V_2$. The dimensions of $V$ and $W$ are denoted by $p$ and
$q$, respectively. One parameter subgroups of $\SL(V)$ will denoted
by the letter $\kappa$, weight vectors used for defining
a one parameter subgroup of $\SL(V)$
will be denoted by $\ul{\delta}=(\delta_1,...,\delta_p)$. 
For one parameter subgroups of $\SL(W)$ we use $\la$,
and write weight vectors as $\ul{\gamma}=(\gamma_1,...,\gamma_q)$.
If we are given bases $v_1,...,v_p$ for $V$ and $w_1,...,w_q$ for $W$,
$i\in\{\, 0,...,p\,\}$ and $j\in\{\,0,...,q\,\}$, we set
$V^{(i)}:=\langle\, 1,...,i\,\rangle$ and $W^{(j)}:=\langle\,1,...,j\,\rangle$.
\par
Now,
 let $(\kappa,\la)$
be an arbitrary one parameter subgroup of $\SL(V)\times\SL(W)$.
Choose bases $v_1,...,v_p$ of $V$ and $w_1,...,w_q$ of $W$ with respect to 
which $\kappa$ and $\la$ act diagonally and are determined
by weight vectors $(\delta_1,...,\delta_p)$ with
$\delta_1\le\cdots\le\delta_p$, $\sum\delta_i=0$,
and $(\gamma_1,...,\gamma_q)$ with $\gamma_1\le\cdots\le\gamma_q$,
$\sum\gamma_j=0$, respectively.
Set 
$j_0 := \min\{\, j\,|\, \Im f\subset W^{(j)}\,\}$,
$i_0 := \min\{\, i\,|\,  f(V^{(i)})\not\subset W^{(j_0-1)}\,\}$,
$i_0^\p := \min \{\,i\,|\, V^{(i)}\not\subset\ker f\,\}$, and
$j_0^\p  := \min\{\,j\,|\,  f(V^{(i_0^\p)})\subset W^{(j)}\,\}$.
Write $ f=\sum f_{i,j} v_i^\vee\otimes w_j$.
Let $s,t$ be indices such that
$f_{s,t}\neq
0$. For $i=1,...,p$, $m(i;s,t)$ denotes the weight of the
eigenvector $v_{s}^\vee \otimes w_{t}$ with respect to the action of
the one parameter
subgroup $\kappa^{(i)}$. In the same way, the numbers $n(j;s,t)$ are
defined. Then, one easily checks:
\begin{Lem}
\label{Obs2}
{\rm i)}
$\mu([f],(\kappa^{(i)},\la^{(0)}))= m(i; s,t)$ unless $i_0^\p\le i < s$.
In that case, we will have $\mu([f],(\kappa^{(i)},\la^{(0)}))=p-i$ and
$m(i;s,t)=-i$.
\par
{\rm ii)}
$\mu([f],(\kappa^{(0)},\la^{(j)}))= n(j; s,t)$ unless $t \le j < j_0$.
Then, $\mu([f],(\kappa^{(0)},\la^{(j)}))= j$ and $n(j; s,t)=j-q$.
\end{Lem}
Theorem~\ref{Decomp1} can now be restated as
\begin{Thm*}
\label{IGIT2}
In the above situation, there exist indices
$i_*$ and
$j_*$ with $f_{i_*,j_*}\neq 0$ and a decomposition  of the weight
vector
\begin{eqnarray*}
(\underline{\delta},\underline{\gamma})
\q =\q 
\sum_{i=1}^{i_0^\p-1} \alpha_i \delta^{(i)} 
&+& 
\sum_{i=i_*}^{p-1}\widetilde{\alpha}_i\delta^{(i)}
\q +\q
\sum_{j=j_0}^{q-1}\beta_j\gamma^{(j)}\q+\\
&+&
\sum_{j=1}^{j_*-1}\widetilde{\beta}_j\gamma^{(j)}
\q+\q\sum_{{i=i_0^\p,...,p-1;\atop
j=1,...,j_0-1\phantom{;}}}
\eta_{i,j} \Bigl({1\over p}\delta^{(i)},{1\over  q}\gamma^{(j)}\Bigr).
\end{eqnarray*}
The $\alpha_i$, $\widetilde{\alpha}_i$, $\beta_j$,
$\widetilde{\beta}_j$, and $\eta_{i,j} $
are non-negative rational numbers such that $\eta_{i,j}=0$
whenever {\rm (a)} $ f(V^{(i)})\not\subset W^{(j)}$, or {\rm (b)}
$i<i_*$ and $j\ge j_*$, or {\rm (c)}
$i\ge i_*$ and $j<j_*$.
\end{Thm*}
\begin{Lem}
\label{Triv1}
Let $\underline{\delta}=\sum_{i=i_1}^{i_2}\alpha_i\delta^{(i)}$,
$\alpha_i\in\Q_{\ge 0}$, $i=i_1,...,i_2$, and
$\underline{\gamma}=\allowbreak\sum_{j=j_1}^{j_2}\beta_j \gamma^{(j)}$,
$\beta_j\in\Q_{\ge 0}$, $j=j_1,...,j_2$, be vectors of rational numbers.
\par
{\rm i)} If $\sum_{i=i_1}^{i_2}p\alpha_i\ge
\sum_{j=j_1}^{j_2}q\beta_j$, then we can decompose the vector
$(\underline{\delta},\underline{\gamma})$
in the following way
$$
(\underline{\delta},\underline{\gamma})\q =\q
\sum_{i=i_1}^{i_2}{\alpha}^\p_i\delta^{(i)}+
\sum_{{i=i_1,...,i_2;\atop
       j=j_1,...,j_2\phantom{;}}} \eta_{i,j}
\Bigl({1\over p}\delta^{(i)},{1\over q}\gamma^{(j)}\Bigr)
$$
where all the ${\alpha}^\p_i$ and $\eta_{i,j}$ are non-negative
rational numbers, $\alpha_i^\p\le\alpha_i$, $i=i_1,...,i_2$, and
$\eta_{i,j}=0$ when
$\alpha_i\cdot\beta_j=0$.
Moreover, if equality holds in the assumption, all the
${\alpha}^\p_i$ are zero.
\par
{\rm ii)} For $\sum_{i=i_1}^{i_2}p\alpha_i\le
\sum_{j=j_1}^{j_2}q\beta_j$, we can write
$$
(\underline{\delta},\underline{\gamma})\q =\q
\sum_{j=j_1}^{j_2}{\beta}^\p_j\gamma^{(j)}+
\sum_{{i=i_1,...,i_2;\atop
       j=j_1,...,j_2\phantom{;}}} \eta_{i,j}
\Bigl({1\over p}\delta^{(i)},{1\over q}\gamma^{(j)}\Bigr)
$$
where all the ${\beta}^\p_j$ and $\eta_{i,j}$ are non-negative
rational numbers, $\beta^\p_j\le \beta_j$, $j=j_1,...,j_2$, and
$\eta_{i,j}=0$ when
$\alpha_i\cdot\beta_j=0$. Furthermore, if equality is assumed, the
${\beta}_j^\p$ are zero.
\end{Lem}
\it Proof\rm. Obvious.
\par
\bigskip
\noindent
Let $\delta_{{\frak i}_1}<\cdots<\delta_{{\frak i}_e}$ 
be the different weights which
appear, and $V:=\bigoplus_{\eps=1}^e V^{\delta_{{\frak i}_\eps}}$
be the corresponding decomposition of $V$ into eigenspaces.
We will write $V^{\le\delta_{{\frak
i}_\eps}}:=\bigoplus_{\eps^\p=1}^\eps
V^{\delta_{{\frak i}_{\eps^\p}}}$ and
$V^{<\delta_{{\frak i}_\eps}}:=\bigoplus_{\eps^\p=1}^{\eps-1}
V^{\delta_{{\frak i}_{\eps^\p}}}$. Set
$\alpha_i:=(1/p)(\delta_{i+1}-\delta_i)$, $i=1,...,p-1$.
\par
Similarly,
let $\gamma_{{\frak j}_1}<\cdots<\gamma_{{\frak j}_f}$
be the different weights occuring. This gives rise to analogous
constructions as before which we will not write down explicitly.
Define $\beta_j:=(1/q)(\gamma_{j+1}-\gamma_j)$, $j=1,...,q-1$.
\par
We have to make the computation at
an eigenvector $v^\vee_{i_*}\otimes w_{j_*}$ with $f_{i_*,j_*}\neq 0$ at
which $\mu([f],(\kappa,\la))$ is achieved. So, we first look at
all weights which appear as weights of eigenvectors appearing
non-trivially in the decomposition of $ f$. Hence, we
define
\begin{eqnarray*}
G:=\Bigl\{\, (\eps,\sigma)\in\{\, 1,...,e\,\}\times\{\,
1,...,f\,\}&|& \exists\ \hat{\imath}\in\{\,1,...,p\,\},\
\hat{j}\in\{\,1,...,q\}:\\
 & & f_{\hat{\imath},\hat{j}}\neq 0\ \wedge\
v_{\hat{\imath}}\in V^{\delta_{{\frak i}_{\eps}}}\ \wedge\
w_{\hat{j}}\in W^{\gamma_{{\frak j}_{\sigma}}}\,\Bigr\}.
\end{eqnarray*}
Let $\gamma_{{\frak k}_1}<\cdots<\gamma_{{\frak k}_{l_0}}$
be the weights with ${\frak k}_\iota={\frak j}_{\sigma}$
for some $(\eps,\sigma)\in G$, $\iota=1,...,l_0$.
Set
${\frak h}_{l_0}:=\min\{\, {\frak i}_\eps\,|\, (\eps,\sigma)\in
G\hbox{ with } \sigma\hbox{ such that } 
{\frak k}_{l_0}={\frak j}_{\sigma}\,\}$.
If $ f(V^{<\delta_{{\frak h}_{l_0}}})=\{0\}$, we stop.
Otherwise, we define $l_1$ by
$$
{\mathfrak k}_{l_1}=
\max\bigl\{\, {\mathfrak j}_\sigma\,|\, (\eps,\sigma)\in G
\hbox{ with } \eps \hbox{ such that } {\mathfrak i}_{\eps+1}\le 
{\mathfrak h}_{l_0}\,\bigr\},
$$
i.e., ${\mathfrak k}_{l_1}$ is determined by the requirements
$
\varphi\bigl(V^{< \delta_{{\mathfrak h}_{l_0}}}\bigr)\subset
W^{\le \gamma_{{\mathfrak j}_\sigma}}
$
and that there be an $(\eps,\sigma)\in G$ with 
${\mathfrak j}_{\sigma}={\mathfrak k}_{l_1}$ and ${\mathfrak i}_{\eps}<{\mathfrak h}_{l_0}$.
Note that $l_1<l_0$. Next, set
${\mathfrak h}_{l_1}:=\min\{\, {\mathfrak i}_\eps\,|\, (\eps,\sigma)\in
G\hbox{ with } \sigma\allowbreak\hbox{ such that } 
{\mathfrak k}_{l_1}={\mathfrak j}_{\sigma}\,\}$.
Now, iterate this process to get --- after relabelling --- weights
$\delta_{{\frak h}_1}<\cdots <\delta_{{\frak h}_l}$
and $\gamma_{{\frak k}_1}<\cdots<\gamma_{{\frak k}_l}$ with the property that
\begin{equation}
\label{eq3}
 f
\bigl(V^{<\delta_{{\frak h}_{1}}}\bigr)\ =\ \{0\},\ \ 
 f\bigl(V^{<\delta_{{\frak h}_{\iota}}}\bigr)\ \subset\  W^{\le
\gamma_{{\frak k}_{\iota-1}}},\ \iota=2,...,l,\ \ \hbox{and}\
 f(V)\ \subset\
W^{\le \gamma_{{\frak k}_{l}}}.
\end{equation}
Finally, we proceed to the proof.
We choose $*\in\{\, 1,...,l\,\}$ such that
$-\delta_{{\frak h}_*}+\gamma_{{\frak k}_*}$ becomes maximal,
so that for any $\iota\in\{\, 1,...,l\,\}$
\begin{equation}
\label{eq4}
-\delta_{{\frak h}_{\iota}}+\gamma_{{\frak k}_\iota}\q\le\q
-\delta_{{\frak h}_{*}}+\gamma_{{\frak k}_*}.
\end{equation}
Let $i^\p$ be minimal with $v_{i^\p}\in V^{\delta_{{\frak h}_*}}$
and $j^\p$ be minimal with $w_{j^\p}\in W^{\gamma_{{\frak k}_*}}$.
Next, write $(\underline{\delta},\underline{\gamma})=
(\underline{\delta}_1,\underline{\gamma}_1)+
(\underline{\delta}_2,\underline{\gamma}_2)$.
Here, $\underline{\delta}_1=\sum_{i=1}^{i^\p-1} \alpha_i\delta^{(i)}$,
$\underline{\gamma}_1=\sum_{j=1}^{j^\p-1}\beta_j \gamma^{(j)}$,
$\underline{\delta}_2=\sum_{i=i^\p}^{p-1} \alpha_i\delta^{(i)}$, and
$\underline{\gamma}_2=\sum_{j=j^\p}^{q-1}\beta_j \gamma^{(j)}$.
We begin by decomposing $(\underline{\delta}_2,\underline{\gamma}_2)$.
First, let $i^{\p\p}$ be maximal with $v_{i^{\p\p}}\in V^{<\delta_{{\frak
h}_{(*+1)}}}$ and $j^{\p\p}$ maximal with $w_{j^{\p\p}}\in
W^{<\gamma_{{\frak k}_{(*+1)}}}$.
Let $\underline{\delta}_2^\p:=\sum_{i=i^\p}^{i^{\p\p}}
\alpha_i\delta^{(i)}$
and $\underline{\gamma}_2^{\p}:=\sum_{j=j^\p}^{j^{\p\p}}
\beta_j\gamma^{(j)}$.
Observe that, by~(\ref{eq3}),  $ f(V^{(i)})\subset W^{(j)}$ for all
$i,j$ with $\alpha_i\cdot\beta_j\neq 0$ (the first non-zero
coeffient of a $\gamma^{(j)}$ is $\beta_{\tilde{j}}$ where
$\tilde{j}$ is maximal such that $w_{\tilde{j}}\in
W^{\le\gamma_{{\frak k}_{*}}}$,
by~(\ref{Decomp0})). By Equation~(\ref{Decomp0}) and~(\ref{eq4})
$$
\sum_{i=i^\p}^{i^{\p\p}} p\alpha_i\q =\q \delta_{{\frak h}_{(*+1)}}
- \delta_{{\frak h}_{*}}\q \ge\q \gamma_{{\frak k}_{(*+1)}}-
\gamma_{{\frak k}_{*}}\q =\q \sum_{j=j^\p}^{j^{\p\p}}q\beta_j.
$$
Thus, we can decompose
$$
(\underline{\delta}_2^\p,\underline{\gamma}_2^\p)\q=\q
\sum_{i=i^\p}^{i^{\p\p}}{\alpha}^\p_i\delta^{(i)}+
\sum_{i=i^\p,...,i^{\p\p};\atop
      j=j^\p,...,j^{\p\p}\phantom{;}}
\eta_{i,j}\Bigl({1\over p}\delta^{(i)},{1\over q}\gamma^{(j)}\Bigl)
$$
according to
Lemma~\ref{Triv1}, i). Moreover, $\eta_{i,j}\neq 0$ implies
$ f(V^{(i)})\subset W^{(j)}$.
Next, define $i^{\p\p\p}$ as the maximal index with
$v_{i^{\p\p\p}}\in V^{<\delta_{{\frak h}_{(*+2)}}}$ and let $j^{\p\p\p}$
be maximal with
$w_{j^{\p\p\p}}\in W^{<\gamma_{{\frak k}_{(*+2)}}}$.
Set $\underline{\delta}^{\p\p}_2:=\sum_{i=i^\p}^{i^{\p\p}}
\alpha_i^\p\delta^{(i)}+\sum_{i=i^{\p\p}+1}^{i^{\p\p\p}}\alpha_i\delta^{(i)}$,
and
$\underline{\gamma}_2^{\p\p}=
\sum_{j=j^{\p\p}+1}^{j^{\p\p\p}}\beta_j\gamma^{(j)}$.
Since, for the same reason as before,
$\sum_{i=i^\p}^{i^{\p\p\p}}p\alpha_i\ge\sum_{j=j^\p}^{j^{\p\p\p}}q\beta_j$,
it is clear that
$
\sum_{i=i^\p}^{i^{\p\p}}
p\alpha_i^\p\allowbreak+\sum_{i=i^{\p\p}+1}^{i^{\p\p\p}}p\alpha_i \ge
\sum_{j=j^{\p\p}+1}^{j^{\p\p\p}}q\beta_j.
$
Again, $ f(V^{(i)})\subset W^{(j)}$ for every $i,j$, such that either
$\alpha^\p_i\beta_j\neq 0$ or $\alpha_i\beta_j\neq 0$.
Hence, we can apply Lemma~\ref{Triv1} again.
\par
Now, iterate this process until all the $\beta_j$'s of the
beginning with $j<j_0$ are eaten up. The result is a decomposition
$$
(\underline{\delta}_2,\underline{\gamma}_2)\q =\q
\sum_{i=i^\p}^{p-1}\widetilde{\alpha}_i\delta^{(i)}
+\sum_{i=i^\p,...,p-1;\atop
       j=j^\p,...,j_0-1\phantom{;}} \eta_{i,j}
\Bigl({1\over p}\delta^{(i)},{1\over q}\gamma^{(j)}\Bigr)+
\sum_{j=j_0}^{q-1}\beta_j\gamma^{(j)},
$$
such that $\eta_{i,j}\neq 0$ implies $ f(V^{(i)})\subset W^{(j)}$.
The corresponding
decomposition of the remaining vector
$$
(\underline{\delta}_1,\underline{\gamma}_1)\q =\q
\sum_{i=1}^{i_0^\p-1} \alpha_i\delta^{(i)}
+\sum_{i=i_0^\p,...,i^\p-1;\atop
       j=1,...,j^\p-1\phantom{;}} \eta_{i,j}
\Bigl({1\over p}\delta^{(i)},{1\over q}\gamma^{(j)}\Bigr)+
\sum_{j=1}^{j^\p-1}\widetilde{\beta}_j\gamma^{(j)}\phantom{,}
$$
is achieved by an analogous method, this time "working backwards"
and making use of~(\ref{eq4}) and Lemma~\ref{Triv1}, ii).
By definition of $G$, we find  indices $i_*$ and $j_*$ such that
$v_{i_*}\in V^{\delta_{{\frak h}_*}}$, $w_{j_*}\in W^{\gamma_{{\frak
k}_*}}$, and $f_{i_*,j_*}\neq 0$. This completes the proof.
\subsubsection*{The general case}
Let
$([w_i], i\in V; [f_a], a\in A)\in \P_{\ul{\tau},\catp}$ and
$\la=(\la_1,...,\la_n)$ be a one parameter subgroup of $\SL(\catP)$.
Choose bases $v_1^i,...,v_{p_i}^i$ for $V_i$, $i=1,...,n$, such that
$\la$ is given w.r.t.\ those bases by the weight vector
$\ul{\gamma}$, and
let $Q(\la)$ be the full subquiver associated with the set
$V(\la):=\{\, i\in V\,|\, \la_i\neq \la_i^{(0)}\,\}$.
Using Theorem~\ref{Decomp1}, for each arrow $a\in A(\la)$,
we pick $f^a_{i(t(a))_*i(h(a))_*}\neq 0$ (in the corresponding
decomposition of $f_a$) and a decomposition of the weight vector
$
(\ul{\gamma}^{t(a)}, \ul{\gamma}^{h(a)})$ as
$$
\sum_{{i(t(a))=0,...,p_{t(a)}-1; \atop
       i(h(a))=0,...,p_{h(a)}-1\phantom{;}}} \eta^0_{i(t(a))i(h(a))}
\Bigl({1\over p_{t(a)}}
\gamma^{t(a),(i(t(a)))}, {1\over p_{h(a)}} \gamma^{h(a),(i(h(a)))}\Bigr)
$$
with the properties asserted in Theorem~\ref{Decomp1}.
Next, suppose we are given a decomposition
$\ul{\gamma}=\sum_{\ul{j}} \eta_{\ul{j}} \gamma^{\ul{j}}$.
This yields, for every arrow $a\in Q(\la)$, a decomposition
of $(\ul{\gamma}^{t(a)}, \ul{\gamma}^{h(a)})$
as
$$
\sum_{{i(t(a))=0,...,p_{t(a)}-1; \atop
       i(h(a))=0,...,p_{h(a)}-1\phantom{;}}} \eta_{i(t(a))i(h(a))}
\Bigl({1\over p_{t(a)}}
\gamma^{t(a),(i(t(a)))}, {1\over p_{h(a)}} \gamma^{h(a),(i(h(a)))}\Bigr)
$$
with
$$
\eta_{i(t(a))i(h(a))}:=\sum_{(j_1,...,j_n):
                                              j_{t(a)}=i(t(a));
                                           \atop\phantom{(j_1,...,j_n):}
                                              j_{h(a)}=i(h(a))\phantom{;}}
\eta_{\ul{j}}.
$$
\begin{Thm}
\label{Decomp2}
There exists a decomposition 
$\ul{\gamma}=\sum_{\ul{j}} \eta_{\ul{j}} \gamma^{\ul{j}}$
into basic weight vectors
with the property that, for every arrow $a\in A$, $i(t(a))=0,...,
p_{t(a)}-1$, and $i(h(a))=0,...,p_{h(a)}-1$   
$$
        \sum_{(j_1,...,j_n):j_{t(a)}=i(t(a));
\atop\phantom{(j_1,...,j_n):}j_{h(a)}=i(h(a))\phantom{;}}
\eta_{\ul{j}}\qquad=\qquad \eta^0_{i(t(a))i(h(a))}.
$$
\end{Thm}
\it Proof\rm .
This will be done by induction, the case
$\#V=2$ being already settled.
In a suitable labelling of the vertices, $\St_Q(n)=\{\, n-1,n\,\}$.
Let us assume for simplicity that the arrow in $\St_Q(n)$ is $(n-1,n)$.
By hypothesis, there exist decompositions
$(\ul{\gamma}^1,...,\ul{\gamma}^{n-1})=\sum_{\ul{j}=(j_1,..,j_{n-1},0)}
\widetilde{\eta}_{\ul{j}} \ul{\gamma}^{\ul{j}}$ and 
$(0,...,0,\ul{\gamma}^{n-1}, \ul{\gamma}^n)=
\sum \eta_{i(n-1)i(n)}^0 \gamma^{(0,...,0,i(n-1),i(n))}$
with the 
respective properties.
Since, by (\ref{Decomp0}), for each $i(n-1)=1,...,p_{n-1}-1$,
$$
\sum_{i(n)} \eta^0_{i(n-1)i(n)}\q=\q
\sum_{(j_1,...,j_{n-1},0):\atop j_{n-1}=i(n-1)} \widetilde{\eta}_{\ul{j}},
$$
the assertion is obvious. 
\par
\bigskip
\noindent
\begin{Cor}
\label{Ess} 
A point $([w_i], i\in V; [f_a], a\in A)$ in $\P_{\ul{\tau},\catp}$
is (semi)stable w.r.t.\ chosen the linearization in $\O(l_1,...,l_n; b_a, a\in A)$,
if and only if the Hilbert-Mumford criterion is fulfilled for all basic one
parameter subgroups.
\end{Cor}
To conclude, let us look at the weigths of the basic one parameter
subgroups. In the above notation, let
$\ul{\gamma}^{\ul{j}}$ be a basic weight vector, and $Q^\p:=Q(\la^{\ul{j}})$.
Then,
\begin{eqnarray*}
\mu\Bigl(([f_a], a\in A), \la^{\ul{j}}\Bigr) &=&
\sum_{a\in A^\p} b_a\Bigl(-{j_{t(a)}\over p_{t(a)}}+{j_{h(a)}\over p_{h(a)}}
\Bigr)
\\
&& 
+\sum_{i\in\ED_Q(Q^\p); a\in \Out_Q(i)} b_a\Bigl(-{j_{t(a)}\over p_{t(a)}}
+\eps_a(V_{t(a)}^{(j_{t(a)})})\Bigr)\\
&&
+\sum_{i\in\ED_Q(Q^\p); a\in \In_Q(i)} b_a\Bigl({j_{h(a)}\over p_{h(a)}}
-\eps_a(V_{h(a)}^{(j_{h(a)})})\Bigr).
\end{eqnarray*}
Here, $\eps_a(V_{t(a)}^{(j_{t(a)})})=1$ if $V_{t(a)}^{(j_{t(a)})}
\not\subset\ker f_a$ and $0$ otherwise, and 
$\eps_a(V_{h(a)}^{(j_{h(a)})})=1$ if $V_{h(a)}^{(j_{h(a)})}
\supset\Im f_a$ and $0$ otherwise.
One easily sees:
\begin{Cor}
\label{Ess2}
Given bases $v^i_1,...,v_{p_i}^i$ for $V_i$, $i\in V$,
and indices
$j_i\in \{\, 0,...,p_i\,\}$ with
$$
f_a\bigl(V_{t(a)}^{(j_{t(a)})}\bigr)\q\subset\q
   V_{h(a)}^{(j_{h(a)})},\qquad\hbox{for all } a\in A,
$$ 
one has
$$
\mu\Bigl(([f_a], a\in A), \la^{(j_1,...,j_n)}\Bigr)\q\le\q
\sum_{a\in A} b_a\Bigl(-{j_{t(a)}\over p_{t(a)}}+{j_{h(a)}\over p_{h(a)}}
\Bigr)
$$
with equality in case $\ul{\gamma}^{(j_1,...,j_n)}$ is basic w.r.t.\
$([f_a], a\in A)$.
\end{Cor}
\subsection*{Proof of Theorem~\ref{GITII}}
If we assume that $l_i=\sum_{a\in A(i)} b_a \alpha_{i, a}$
for some positive rational numbers $\alpha_{i,a}$, we find in the
setting of Corollary~\ref{Ess2}
\begin{eqnarray*}
&&\mu\Bigl(\bigl([w_i], i\in V;[f_a], a\in A\bigl), \la^{{(j_1,...,j_n)}}\Bigr) \\
&\le&\sum_{a\in A} b_a\Bigl[\alpha_{t(a),a}{1\over p_{t(a)}}
\mu\bigl([w_{t(a)}], \la_{t(a)}^{(j_{t(a)})}\bigl)-{j_{t(a)}\over p_{t(a)}}
\\
&&
\phantom{\sum_{a\in A} b_a\Bigl[
} + \alpha_{h(a),a}{1\over p_{h(a)}}
\mu\bigl([w_{h(a)}], \la_{h(a)}^{(j_{h(a)})}\bigr)+{j_{h(a)}\over p_{h(a)}}
\Bigr]
\end{eqnarray*}
with equality in case $\ul{\gamma}^{(j_1,...,j_n)}$ is basic w.r.t.\
$([f_a], a\in A)$.
The assertion now follows from 
Corollary~\ref{Ess}.
\begin{Rem}
Observe that we made use of the Additivity Assumption~\ref{Add}
in this computation.
It is not valid without it.
\end{Rem}
\section{Proof of Theorem~\ref{MTH}}
We can now start with the GIT construction which follows a well-known
pattern. For the case of semistable torsion free coherent sheaves, one
may consult \cite{HL}. The details left out here can be easily filled
in with that reference.
\subsection*{Boundedness}
\label{bound}
We must first convince ourselves that the semistable objects live
in bounded families, i.e., can be parametrized by a scheme of finite
type over $\C$. Here, we check that this is true for the participating
torsion free coherent sheaves.
\par
The data $\catP$, $\ul{\sigma}_Q$, and $\ul{b}_Q$ are the same  
as before, and $\vartheta:=\vartheta(\catP,\ul{\sigma}_Q,\ul{b}_Q)$. 
Let $\catR=(\E_i,i\in V; \phi_a,a\in A)$ be a
$\vartheta$-semistable representation
of type $\catP$, and $i_0\in V$ a fixed vertex. 
Observe that $\Gamma_Q$'s being a tree implies that every other vertex can
be connected to $i_0$ by a unique (shortest) path.  
For any non-trivial, proper subsheaf
$\G$ of $\E_{i_0}$, we define the  sub-representation 
$\catR_\G=(\F_i, i\in V;\phi_a^\p, a\in A)$ of $\catR$ as follows: 
We set $\F_{i_0}:=\G$, and for any other
vertex $i$, we define $\F_i$ as $0$ if the path connecting $i$ to $i_0$
passes through an ingoing arrow of $i_0$ and as $\E_i$ in the other case,
and finally, $\phi_a^\p:=\phi_{a|\F_{t(a)}}$ for all $a\in A$.
The condition of $\vartheta$-semistability applied
to $\catR_\G$ shows that
\begin{eqnarray*}
&\sum_{a\in \mathop{\rm Out}_Q(i_0)}&
b_a\biggl[ \check{\sigma}_{i_0}\Bigl\{ P(\F_{i_0})-\rk\F_{i_0}
\Bigl({P(\E_{i_0})-\sigma_{i_0}\over\rk\E_{i_0}}+{\sigma_{i_0}\over\rk\F_{i_0}}
\Bigr)\Bigr\}\biggr]\\
&+\sum_{a\in \mathop{\rm In}_Q(i_0)}&
b_a\biggl[ \check{\sigma}_{i_0}\Bigl\{ P(\F_{i_0})-\rk\F_{i_0}
{P(\E_{i_0})+\sigma_{i_0}\over\rk\E_{i_0}}
\Bigr\}\biggr]
\end{eqnarray*}
is a non-positive polynomial, whence
\begin{equation}
\label{boundx}
\Bigl(\sum_{a\in A(i_0)} b_a\check{\sigma}_{i_0}\Bigr)
{P(\F_{i_0})\over\rk\F_{i_0}}\q\le\q 
\Bigl(\sum_{a\in A(i_0)} b_a\check{\sigma}_{i_0}\Bigr)
{P_{i_0}\over r_{i_0}} +
{\sum_{a\in A(i_0)} b_a{\sigma}}.
\end{equation}
Let $\ol{\sigma}_{i_0}^{{\vee}}$ be the leading coefficient of
the polynomial $\sum_{a\in A(i_0)} b_a\check{\sigma}_{i_0}$,
$\check{s}_{i_0}$ its degree,
$\ol{\sigma}_{i_0}$ the coefficient of
the term of degree $(\check{s}_{i_0}+\dim X-1)$ in $(\sum_{a\in A}b_a)\sigma$
(the latter polynomial has degree at most $(\check{s}_{i_0}+\dim X-1)$),
and define the constant
$C_{i_0}:=\ol{\sigma}_{i_0}/\ol{\sigma}_{i_0}^{{\vee}}$.
Thus, taking leading coefficients in (\ref{boundx}) shows
$$
\mu_{\max}(\E_{i_0})\q\le\q \mu_{i_0}+C_{i_0},
$$
whence, by invoking the boundedness theorem of Maruyama \cite{HL}, 
Theorem~3.3.7,
the following result is obtained:
\begin{Thm}
\label{bound1}
For every $i\in V$, the set of isomorphy classes of torsion free
coherent sheaves $\E_i$ with Hilbert polynomial $\catP(i)$ 
showing up in $\vartheta$-semi\-stable 
representations
of type $\catP$ is bounded.
\end{Thm}
\subsection*{Sectional semistability}
Sectional semistability is a technical way of rewriting the semistability
condition. It is the form in which we will encounter it
during the GIT construction. 
\par
This time, we fix in addition to $\catP$ and $\ul{b}_Q$
positive rational numbers $s_i$, $i\in V$,
and set $s:=s_1\cdot...\cdot s_n$, $\check{s}_i=s/s_i$, $i\in V$.
A representation $(\E_i, i\in V;\phi_a, a\in A)$ of $Q$ is called
\it $(\ul{s}_Q,\ul{b}_Q)$-sectional (semi)stable\rm , 
if there are subspaces $H_i\subset H^0(\E_i)$ of dimension $\chi(\E_i)$,
$i\in V$, 
with $f_a(H_{t(a)})\subset H_{h(a)}$, for all $a\in A$, such that,
for all sub-representations $\ul{\F}=(\F_i,i\in V, \phi_a^\p, a\in A)$, the
number
\begin{eqnarray*}
\delta(\ul{\F},\ul{\E},\ul{\phi}) &:=&
\sum_{a\in A} b_a\biggl[\check{s}_{t(a)}\Bigl\{
\dim \bigl(H_{t(a)}\cap H^0(\F_{t(a)})\bigr)-
\rk\F_{t(a)}{\chi(\E_{t(a)})-s_{t(a)}\over\rk\E_{t(a)}}
\Bigr\}\\
& + &
\check{s}_{h(a)}\Bigl\{
\dim\bigl(H_{h(a)}\cap H^0(\F_{h(a)})\bigr)
-\rk\F_{h(a)}{\chi(\E_{h(a)})+s_{h(a)}\over\rk\E_{h(a)}}
\Bigr\}\biggr]
\end{eqnarray*}
is negative (non-positive).
\begin{Thm}
\label{bound10}
There exists an $m_0$, such that, 
for all $i_0\in V$, the set of isomorphy classes of torsion
free coherent sheaves $\E_{i_0}$ occuring in representations
$(\E_i,i\in V; \phi_a,a\in A)$ of type $\catP$
with the property that there exists an $m\ge m_0$
such that
$(\E_i(m), i\in V;\phi_a(m), a\in A)$ is sectional semistable 
w.r.t.\ the parameters $(\sigma_i(m), i\in V; b_a, a\in A)$
is bounded, too.
\end{Thm}
{\it Proof}.
Before starting with the proof, we remind the reader
of the following important result.
\begin{Thm*}[\rm The Le Potier-Simpson estimate \cite{HL}, p.71]
Let $\F$ be a torsion free coherent sheaf, 
and $C(\F):=\rk\F(\rk\F+\dim X)/2$.
Then, for every $m\ge 0$,
\begin{eqnarray*}
{h^0(\F(m))\over\rk \F}&\le& {\rk\F-1\over(\dim X)!
\rk\F}\bigl[\mu_{\max}(\F)+C(\F)-1+m\bigr]_+^{\dim X}\\
& &+\ {1\over(\dim X)!\rk\F}
\bigl[\mu(\F)+C(\F)-1+m\bigr]_+^{\dim X}.
\end{eqnarray*}
\end{Thm*}
As in the proof of Theorem~\ref{bound1}, 
for $i_0\in V$ and a non-trivial subsheaf $\G$
of $\E_{i_0}$, 
we find, for $m\gg 0$, an upper bound
$$
{h^0(\G(m))\over \rk\G}\q\le\q {\chi(\E_{i_0}(m))\over\rk \E_{i_0}}
+K_{i_0}(m)\q\le\q
{h^0(\E_{i_0}(m))\over\rk \E_{i_0}}
+K_{i_0}(m).
$$
Here, $K_{i_0}(m)$ is a positive rational function growing at most like 
a polynomial of degree $\dim X-1$ depending only on the input data
$\ul{\sigma}_Q$, $\ul{b}_Q$, and $\catP$.
Hence, for every non-trivial quotient ${\cal Q}$
of $\E_{i_0}$, we get a lower
bound
\begin{eqnarray*}
{\chi(\E_{i_0}(m))\over\rk \E_{i_0}}
-(\rk\E_{i_0}-1)K_{i_0}(m)&\le&
{h^0(\E_{i_0}(m))\over\rk \E_{i_0}}
-(\rk\E_{i_0}-1)K_{i_0}(m)\\
&\le&
{h^0({\cal Q}(m))\over \rk{\cal Q}}.
\end{eqnarray*}
As remarked before, the left hand side is positive for all
$m\gg 0$ and depends only on the input data
$\ul{\sigma}_Q$, $\ul{b}_Q$, and $\catP$, and we can bound $K_{i_0}$
from above by a polynomial $k_{i_0}(m)$ of degree at most $\dim X-1$.
Thus, applying the above estimate to the minimal slope destabilizing
quotient sheaf of $\E_{i_0}$ and making use of the Le Potier-Simpson
estimate, we find 
\begin{equation}
\label{secti}
{P_{i_0}(m)\over r_{i_0}}-(r_{i_0}-1)k_{i_0}(m)\q\le\q
{1\over\dim X!}\bigl[\mu_{\min}(\E_{i_0})+C-1+m\bigr]_+^{\dim X},
\end{equation}
$C:=\max\{\, s(s+\dim X)/2\,|\, s=1,...,r_{i_0}-1\,\}$.
Let $D$ be the coefficient of $m^{\dim X-1}$ in the left hand polynomial
multiplied by $(\dim X)!$. Then, we can find an $m_{i_0}$ such that
for all $m\ge m_{i_0}$,
$$
{P_{i_0}(m)\over r_{i_0}}-(r_{i_0}-1)k_{i_0}(m)\q>\q
{1\over\dim X!}\bigl[D-{1\over r!}+m\bigr]_+^{\dim X}.
$$
Thus, from (\ref{secti}) and the assumption that 
there exist an $m\ge m_{i_0}$
such that
$(\E_i(m), i\in V;\allowbreak \phi_a(m), a\in A)$ is sectional semistable 
w.r.t.\ the parameters $(\sigma_i(m), i\in V;\allowbreak b_a, a\in A)$,
one finds the lower bound $\mu_{\min}(\E_{i_0})\ge D-C+1$
and, consequently, an upper bound for $\mu_{\max}(\E_{i_0})$,
so that the theorem follows again from Maruyama's boundedness
theorem.
\par\bigskip\noindent
Our next task will be to describe the relation between sectional
semistability for large $m$ and $\vartheta$-semistability.
For this, we have to remind the reader of the following:
\par
Let $\E$ be a torsion free coherent sheaf. A subsheaf $\F\subset\E$ is 
called \it saturated\rm, if $\E/\F$ is again torsion free.
The sheaf $\widetilde{\F}:=\ker\bigl(\E\lra (\E/\F)/\mathop{\rm Tors}(\E/\F)
\bigr)$ is called the \it saturation of $\F$\rm.
It has the same rank as $\F$ and $P(\F)\le P(\widetilde{\F})$.
It is often necessary to restrict to saturated subsheaves, because
of the following boundedness result:
\begin{Thm*}
Let ${\frak C}$ be a bounded family of torsion free coherent sheaves
on $X$. Given a constant $C$, the set of isomorphism classes
of sheaves $\F$ for which there exist a sheaf $\E$ with $[\E]\in {\frak C}$
and a saturated subsheaf $\F^\p$ of $\E$ with $C\le \mu(\F^\p)$
and $\F\cong \F^\p$ is also bounded.
\end{Thm*}
\it Proof\rm.
Note that one can find a finite dimensional vector space $U$ and
an integer $m$, such that every sheaf $\E$ with $[\E]\in {\frak C}$
can be written as a quotient of $U\otimes\O_X(-m)$.
Our assumption yields a lower bound $\mu(\E/\F^\p)\ge C^\p$
for some $C^\p$ which can be computed from $C$ and the maximal slope
occuring for a sheaf in ${\frak C}$.
Grothendieck's lemma (\cite{HL}, Lem.~1.7.9) shows that the family of sheaves
of the form $\E/\F^\p$ (viewed as a quotient of $U\otimes\O_X(-m)$)
is bounded, whence also the family of sheaves of the form $\F^\p$.
\par
\medskip
\noindent  
Now, let $\catR=(\E_i,i\in V;\phi_a,a\in A)$ be a representation of $Q$
and $\catR^\p:=(\F_i,i\in V;\phi^\p_a, a\in A)$ a sub-representation.
Then, for every arrow $a\in A$,
$\phi_a(\widetilde{\F}_{t(a)})\subset \widetilde{\F}_{h(a)}$, so that
$\widetilde{\catR}^\p:=(\widetilde{\F}_{i},i\in V; 
\phi_{a|\widetilde{\F}_{t(a)}}, a\in A)$ is another sub-representation
of $\catR$ with $\vartheta(\catR^\p)\le\vartheta(\widetilde{\catR}^\p)$.
Therefore, $\vartheta$-(semi)stability needs to be checked only for
\it saturated sub-representations $\catR^\p$\rm, i.e.,
sub-representations where $\F_i$ is a saturated subsheaf of $\E_i$,
$i\in V$.  
\begin{Prop}
\label{bound11}
There exists a number $m_1\ge m_0$, such that for every $m\ge m_0$, a   
representation 
$(\E_i,i\in V; \phi_a, a\in A)$
of type
$\catP$
is $\vartheta$-(semi)stable, if and only if 
$(\E_i(m), i\in V;\phi_a(m), a\in A)$ is sectional (semi)stable 
w.r.t.\ the parameters $(\sigma_i(m), i\in V; b_a, a\in A)$ .
\end{Prop}
{\it Proof}.
We show the ``only if''-direction and take
$H_i=H^0(\E_i(m))$. The other direction is similar
but easier, and, hence, omitted.
First, let $\E$ be an arbitary torsion free coherent sheaf of slope $\mu$
and rank $r$,
and suppose we are given a bound $\mu_{\max}(\E)\le \mu+ K$,
for some constant $K$. 
Let $\F$ be a subsheaf with $\mu(\F)\le \mu-rC(\E)-(r-1)K-K^\p$,
for some other constant $K^\p$. Then, the Le Potier-Simpson estimate
tells us that for large $m$
$$
{h^0(\F(m))\over\rk \F}\q\le\q {m^{\dim X}\over {\dim X}!}
+{m^{\dim X-1}\over (\dim X-1)!} (\mu-K-1)+\hbox{lower order terms}.
$$
Now, suppose we are
given some constants $K_i$, $i\in V$ (not to be confounded with
the rational functions of the same name we have considered before).
Then, by the above, for every sub-representation 
$(\F_i,i\in V; \phi_a^\p, a\in A)$
with $\mu(\F_i)\le \mu_i-r_iC(\E_i)-(r_i-1)C_i-K_i$, the $C_i$
being the constants obtained in the proof of Theorem~\ref{bound1}, we get
an estimate
$$
\delta\Bigl(\ul{\F}(m),\ul{\E}(m),\ul{\phi}(m)\Bigr)\q
\le\q
\delta\bigl(\ul{\sigma}_Q,\ul{b}_Q,P_i, K_i, i\in V\bigl)
$$
where the right hand side is a polynomial of degree at most $\dim X-1$
depending only on the data in the bracket. Due to the fact that
the $\sigma_i$, $i\in V$, may have different degrees, that polynomial
might not be homogeneous. First, we define $V_1:=\{\, i\in V\,|\,
\check{\sigma}_i\hbox{ has maximal degree}\,\}$.
Given the fact that we have already uniformly bounded $\mu_{\max}(\E_i)$
for all $i\in V$, it is now clear that we can choose the constants
$K_i$, $i\in V_1$, in such a way 
that for every sub-representation $(\F_i,i\in V;
\phi^\p_a, a\in A)$
with $\mu(\F_{i_0})\le \mu_{i_0}-r_{i_0}C(\E_{i_0})-(r_{i_0}-1) C_{i_0}
-K_{i_0}$ for \sl one single \rm index $i_0\in V_1$, the function
$\delta(\ul{\F}(m),\ul{\E}(m),\ul{\phi}(m))$ will be negative
for all large $m$.
Hence, we may restrict our attention to saturated sub-representations
$(\F_i,i\in V; \phi_a^\p, a\in A)$ where 
$\mu(\F_{i}) > \mu_{i}-r_{i}C(\E_{i})-(r_{i}-1) C_{i}
-K_{i}$ for all $i\in V_1$.
But this means that the $\F_i$, $i\in V_1$, vary in bounded families,
whence they can all be assumed to be globally generated and without 
higher cohomology.
Under these circumstances, the fact that $(\E_i,i\in V;\phi_a, a\in A)$
is $\vartheta$-semistable implies that
the sum of the contributions to $\delta(\ul{\F}(m),\ul{\E}(m),\ul{\phi}(m))$
of the terms coming from $\sigma_i$'s with $i\in V_1$ cannot be positive
for large $m$.
Therefore, we look at $V_2:= \{\, i\in V\,|\,
\check{\sigma}_i\hbox{ has second largest degree}\,\}$.
By iterating the above procedure, we finally get the claim.
\par\bigskip\noindent
The above argumentation shows that the sub-representations of interest
live in bounded families, whence we conclude
\begin{Thm}
\label{bound22}
There is a number $m_2\ge m_1$, such that for all $m\ge m_2$
the following conditions on a representation 
$\catR=(\E_i,i\in V;\phi_a, a\in A)$ of type $\catP$
are equivalent.
\begin{enumerate}
\item $\catR$ is $\vartheta$-(semi)stable.
\item $(\E_i(m),i\in V;\phi_a(m), a\in A)$ fulfills the condition
      $(\sigma_i(m), i\in V;\ul{b}_Q)$-sectional (semi)\allowbreak stabililty
      for all sub-representations $(\F_i(m),i\in V; \phi_a^\p(m), a\in A)$, 
      such that $\F_i(m)$
      is globally generated for all $i\in V$.
\end{enumerate}
\end{Thm}
\subsection*{The parameter space}
By Theorem~\ref{bound1}, we can find a natural number $m_3$
which we choose larger than $m_2$ in Theorem~\ref{bound22},
such that for every $m\ge m_3$ and every 
$\vartheta$-semistable
representation $(\E_i,i\in V; \phi_a,a\in A)$ 
of type $\catP$, the sheaves
$\E_i(m)$ are globally generated and without higher cohomology, $i\in V$.
Fix such an $m$, set $p_i:=P_i(m)$, and choose 
complex vector spaces $V_i$ of dimension $p_i$, $i\in V$.
For every $i\in V$, let ${\frak Q}^\p_i$ be the quasi-projective 
quot scheme of equivalence classes of quotients $k\colon V_i\otimes\O_X(-m)
\lra \E_i$ where $\E_i$ is a torsion free coherent sheaf with Hilbert
polynomial $P_i$ and $H^0(k(m))$ is an isomorphism, and define ${\frak Q}_i$
as the union of those components of ${\frak Q}_i^\p$ which contain
members of $\vartheta$-semistable representations.
\par
For a point $([k_i\colon V_i\otimes\O_X(-m)\lra \E_i, i\in V])\in 
\prod_{i\in V}{\frak Q}_i$ and  an arrow $a\in A$, an element
$\phi_a\in\Hom(\E_{t(a)},\E_{h(a)})$ is uniquely determined by the
associated homomorphism $f_a:=H^0(k_{h(a)}(m))^{-1}\circ H^0(\phi_a(m))\circ
H^0(k_{t(a)}(m))\in \Hom(V_{t(a)}, V_{h(a)})$. 
Define the spaces ${\frak P}:=\prod_{a\in A}\P(\Hom(V_{t(a)}, V_{h(a)})^\vee)$,
${\frak U}:=(\prod_{i\in V} {\frak Q}_i)\times {\frak P}$, 
denote by ${\frak N}_a$ the pullback of the sheaf
$\O_{\P(\Hom(V_{t(a)},V_{h(a)})^\vee)}(1)$ to ${\frak U}\times X$,
and let 
$$
{\frak f}_a\colon V_{t(a)}\otimes\O_{{\frak U}\times X}\lra 
V_{h(a)}\otimes {\frak N}_a
$$
be the pullback of the universal homomorphism on 
$\P(\Hom(V_{t(a)},V_{h(a)})^\vee)$ to 
${\frak U}\times X$.
Moreover, let $\widetilde{\frak k}_i\colon 
V_i\otimes \O_{{\frak U}\times X}\allowbreak
\lra 
\widetilde{\frak E}_i$
be the pullbacks of the universal quotients
on ${\frak Q}_i\times X$ twisted by
$\id_{\pi_X^*\O_X(m)}$, $i\in V$.
Composing ${\frak f}_a$ and $\widetilde{\frak k}_{h(a)}
\otimes {\frak N}_a$, 
we obtain a homomorphism
$$
\widetilde{\frak f}_a\colon
V_{t(a)}\otimes \O_{{\frak U}\times X}\lra 
\widetilde{\frak E}_{h(a)}\otimes {\frak N}_a.
$$
We define ${\frak T}$ as the closed subscheme whose closed points
are those $u\in {\frak U}$ for which $\widetilde{\frak f}_{a|\{u\}\times X}$
vanishes on $\ker(\widetilde{\frak k}_{t(a)|\{u\}\times X})$ for all
$a\in A$, i.e., those
$u$ for which the restriction
$\widetilde{\frak f}_{a|\{u\}\times X}$ comes from
a homomorphism $\phi_{u,a}\colon \widetilde{\frak E}_{t(a)|\{u\}\times X}(-m)
\allowbreak\lra\allowbreak \widetilde{\frak E}_{h(a)|\{u\}\times X}(-m)$.
The scheme ${\frak T}$ is the common vanishing locus of the
vector bundle maps 
$$
\pi_{{\frak U}*}\Bigl(\ker \widetilde{\frak k}_{t(a)}\otimes \pi_{X}^*\O_X(l)\Bigr)
\lra
\pi_{{\frak U}*}\Bigl(\widetilde{\frak E}_{h(a)}\otimes 
{\frak N}_a\otimes \pi_{X}^*\O_X(l)\Bigr),\q a\in A,
$$
obtained
by projecting $\widetilde{\frak f}_a\otimes\pi_X^*\O_X(l)$ 
for $l$ large enough.
Set ${\frak E}_{{\frak T},i}:=
\widetilde{\frak E}_{i|{\frak T}\times X}$, $i\in V$, and let
$\phi_{{\frak T},a}\colon{\frak E}_{{\frak T},t(a)}
\lra {\frak E}_{{\frak T},h(a)}\otimes 
{\frak N}_{a|{\frak T}\times X}$, $a\in A$, be the induced homomorphisms.
The family $({\frak E}_{{\frak T},i},i\in V;
\phi_{{\frak T},a}, a\in A)$ will be abusively called
the \it universal family \rm though it is not exactly a family
in the sense of our definition.
\par
Recall that, for each $i\in V$, there is a left action of
$\SL(V_i)$ on ${\frak Q}_i$ given on closed
points by 
$$
g\cdot [k\colon V_i\otimes\O_X(-m)\lra \E_i]
\q=\q[V_i\otimes \O_X(-m)\stackrel{g^{-1}\otimes{\id}}{\lra}
V_i\otimes \O_X(-m)\stackrel{k}{\lra}\E_i].
$$
Thus, there is an induced left action of $G:=\prod_{i\in V}\SL(V_i)$ on 
$(\prod_{i\in V}{\frak Q}_i)
\times {\frak P}$ which leaves the parameter space
${\frak T}$ invariant and thus yields an action
$$
\alpha\colon G\times {\frak T}\lra {\frak T}.
$$
The proof of the following is standard (cf., e.g., \cite{OST}) 
and, therefore, will not be given 
here.
\begin{Prop}
\label{standard}
{\rm i)} The space ${\frak T}$ enjoys the {\rm local
universal property}, i.e., for every noetherian scheme $S$
and every family $(\ul{\frak E}_S,\ul{\phi}_S)=
({\frak E}_{S, i}, i\in V; \phi_{S,a}, a\in A)$
of $\vartheta$-semistable representations
of type $\catP$ para\-metrized by $S$, there exists
an open covering $U_1,...,U_s$ of $S$ and morphisms 
$h_\nu\colon U_\nu \lra {\frak T}$, such that the
pullback of the universal family on ${\frak T}\times X$
via $(h_\nu\times {\id}_X)$ is equivalent to $(\ul{\frak E}_{S|U_\nu},
\ul{\phi}_{S|U_\nu})$, $\nu=1,...,s$.
\par
{\rm ii)} Given a noetherian scheme $S$ and two morphisms
$h^1$, $h^2$, such that the pullbacks of ${\frak N}_a$, $a\in A$,
via $(h^1\times {\id}_X)$ and $(h^2\times {\id}_X)$ are
trivial and the pullbacks
of the universal family via $(h^1\times {\id}_X)$ and $(h^2\times {\id}_X)$
are equivalent, there exists an \'etale covering
$\eta\colon T\lra S$ and a morphism $t\colon T\lra G$
with $(h^1\circ \eta)=t\cdot (h^2\circ \eta)$.
\end{Prop}
As we will see below, the set ${\frak T}^{\vartheta-(s)s}_{\catp}
\subset {\frak T}$ 
parametrizing $\vartheta$-(semi)stable
representations is open. For this reason and by the universal
property of the good (geometric) quotient, i) and ii) of Theorem~\ref{MTH}
are a direct consequence of the following 
\begin{Thm}
\label{quot10}
The good quotient ${\frak T}^{\vartheta-ss}_{\catp}\catqot G$ 
exists as projective scheme, and the open
subscheme ${\frak T}^{\vartheta-s}_{\catp}\catqot G$
is a geometric quotient for
${\frak T}^{\vartheta-s}_{\catp}$ 
w.r.t.\ the action of $G$.
\end{Thm}
We have chosen to prove this result via Gieseker's method
for the single reason that after determining the weights we are
done and don't have to make any more considerations
about polynomials and leading coefficients which, in our situation,
would be rather messy, I guess.
\subsection*{The Gieseker space and the Gieseker map}
For $i\in V$, let ${\frak k}_i\colon V_i\otimes\pi_X^*\O_X(-m)
\lra {\frak E}_{{\frak Q}_i}$ be the universal quotient
on ${\frak Q}_i\times X$.
The line bundle $\det({\frak E}_{{\frak Q}_i})$ induces a morphism
${\frak d}_{{\frak Q}_i}\colon {\frak Q}_i\lra \Pic X$.
Let ${\frak A}_i$ be the union of the finitely many components of 
$\Pic X$ hit by this map. Note that this does not depend on the
choice of $m$. Therefore, in addition to the other hypothesis' on $m$,
we may assume that $\L(r_im)$ is globally generated
and without higher cohomology, for every $[\L]\in {\frak A}_i$, $i\in V$.
Fix a Poincar\'e line bundle $\n$ on $\Pic X\times X$, let
$\n_i$ be its
restriction to ${\frak A}_i\times X$, and set
$$
{\Bbb G}_i\q:=\q \P\Bigl( \ul{\Hom}\bigl(\bigwedge^{r_i}V_i\otimes
\O_{{\frak A}_i}, \pi_{{\frak A}_i*} 
(\n_i\otimes \pi_X^*\O_X(r_im))\bigr)^\vee\Bigr).
$$
The homomorphism
$$
\pi_{{\frak Q}_i*}\bigl(\mathop{\wedge}^{r_i}(
{\frak k}_i\otimes {\id}_{\pi_X^*\O_X(m)})\bigr)
\colon \bigwedge^{r_i} V_i\otimes\O_{{\frak Q}_i}\lra
\pi_{{\frak Q}_i*}\Bigl(\det({\frak E}_{{\frak Q}_i})\otimes\pi_X^*\O_X(r_im)
\Bigr)
$$
gives rise to an injective and $\SL(V_i)$-equivariant morphism
$\mathop{\rm Gies}_i\colon {\frak Q}_i\lra {\Bbb G}_i$.
Defining ${\frak G}:= (\prod_{i\in V} {\Bbb G}_i)\times {\frak P}$,
we obtain from the above data a $G$-equivariant and injective 
morphism
$$
{\mathop{\rm Gies}}\colon {\frak T}\lra {\frak G}.
$$
We linearize the $G$-action on the right hand space
in $\O(l_1,...,l_n; b_a, a\in A)$ where
$$
l_i\q :=\q \sum_{a\in \Out_Q(i)} b_a{p_{t(a)}-\sigma_{t(a)}(m)
\over r_{t(a)}\sigma_{t(a)}(m)}+\sum_{a\in \In_Q(i)}
b_a{p_{h(a)}+\sigma_{h(a)}(m)\over r_{h(a)}\sigma_{h(a)}(m)}.
$$
After these preparations, Theorem~\ref{quot10} and Theorem~\ref{MTH}
will follow from
\begin{Thm}
{\rm i)} The image of a point
$([k_i\colon V_i\otimes\O_X(-m)\lra \E_i], i\in V; [\phi_a],\allowbreak 
a\in A)$
under $\mathop{\rm Gies}$ in ${\frak G}$
is (semi/poly)stable w.r.t.\ the fixed linearization if and only if
$(\E_i, i\in V_i; \phi_a, a\in A)$ is a 
$\vartheta$-(semi/poly)stable representation of $Q$ of type $\catP$.
\par
{\rm ii)} The morphism
$\mathop{\rm Gies}^{ss}
\colon {\frak T}^{\vartheta-ss}_{\catp}
\lra {\frak G}^{ss}$ is injective and proper and, therefore, finite.
\end{Thm}
\it Proof\rm. i) 
The Gieseker space ${\frak G}$  maps $G$-invariantly
to $\prod_i{\frak A}_i$,
and the fibres are closed, $G$-invariant subschemes, namely,
over the point $\ul{\L}:=(\L_1,...,\L_n)\in \prod_i{\frak A}_i$ sits the
space
$$
{\frak G}_{\ul{\L}}\q=\q
\prod_i\P\Bigl( \Hom\bigl(\bigwedge^{r_i}V_i, H^0(\L_i(r_im))\bigr)^\vee\Bigr)
\times {\frak P}.
$$
It suffices to determine the (semi)stable points in those spaces.
In the following, we use the notation of Section~\ref{GITI}.
We need two formulas.
\begin{itemize}  
\item Let $i\in V$ and $[k_i\colon V_i\otimes\O_X(-m)\lra\E_i]$ be a point
      in ${\frak Q}_i$ with image 
      $g_i:=\mathop{\rm Gies}_i([k_i])\in {\Bbb G}_i$.
      Given a basis $v_1^i,..., v^i_{p_i}$ of $V_i$ and an index 
      $j_i\in\{\, 0,...,p_i\,\}$, we define
      $$
      \E_i^{(j)}\q:=\q k_i\bigl(V_i^{(j)}\otimes \O_X(-m)\bigr).
      $$
      Recall that one has
      \begin{equation}
      \label{weighti}
      \mu\bigl(g_i,\la_i^{(j)}\bigr)\q=\q p_i\rk\E^{(j)}_i-j r_i
      \end{equation}
      and that Assumption~\ref{Add} is verified on the subset
      $\mathop{\rm Gies}_i({\frak Q}_i)\subset {\Bbb G}_i$.
\item Let $a\in A$ be an arrow. Suppose we are given bases
      $v_1^i,...,v_{p_i}^i$ of $V_i$ for $i=t(a)$ and $i=h(a)$.
      Then, one has the following identities
      \begin{equation}
      \label{weightii}
      \begin{array}{ccc}
      {p_{t(a)}-\sigma_{t(a)}(m)\over r_{t(a)}\sigma_{t(a)}(m)}
      {1\over p_{t(a)}}&\bigl(&
       p_{t(a)}\rk \E_{t(a)}^{(j_{t(a)})}-j_{t(a)}r_{t(a)}\q\bigr)\q
      -{j_{t(a)}\over p_{t(a)}}\\ 
       \\
       &=& {p_{t(a)}\rk \E_{t(a)}^{(j_{t(a)})}\over r_{t(a)}
       \sigma_{t(a)}(m)} - {\rk \E_{t(a)}^{(j_{t(a)})}\over r_{t(a)}}
       -{j_{t(a)}\over \sigma_{t(a)}(m)}\\
       \\
       \hbox{and}\hfill
       \\
       \\
       {p_{h(a)}+\sigma_{h(a)}(m)\over r_{h(a)}\sigma_{h(a)}(m)}
       {1\over p_{h(a)}}&\bigl(&
       p_{h(a)}\rk \E_{h(a)}^{(j_{h(a)})}-j_{h(a)}r_{h(a)}\q\bigr)\q
       +{j_{h(a)}\over p_{h(a)}}\\ 
       \\
       &=& {p_{h(a)}\rk \E_{h(a)}^{(j_{h(a)})}\over r_{h(a)}
       \sigma_{h(a)}(m)} + {\rk \E_{h(a)}^{(j_{h(a)})}\over r_{h(a)}}
       -{j_{h(a)}\over \sigma_{h(a)}(m)}.
       \end{array}
       \end{equation}
\end{itemize}
Let $\catt:=\bigl([k_i\colon V_i\otimes\O_X(-m)\lra\E_i],i\in V; [\phi_a],
a\in A\bigr)$ and $\catg:=\mathop{\rm Gies}(\catt)
=\bigl([w_i],i\in V; [f_a], a\in A\bigr)\in {\frak G}$.
\par
By Theorem~\ref{GITII}
and Formula (\ref{weighti}), 
the point $\catg$ is (semi)stable w.r.t.\ the
linearization in $\O(l_1,...,l_n;b_a, a\in A)$, if and only if, for every
possible choice of bases $v_1^i,...,v_{p_i}^i$ for $V_i$ and
indices $j_i\in\{\,0,...,p_i\,\}$
with
$$
f_a\bigl(V_{t(a)}^{(j_{t(a)})}\bigr)\q\subset\q V_{h(a)}^{(j_{h(a)})}
\q\hbox{for every arrow } a\in A,
$$
one has 
\begin{equation}
\label{weightiii}
\begin{array}{c}
0(\ge)\sum_{a}b_a\biggl[
      {p_{t(a)}-\sigma_{t(a)}(m)\over r_{t(a)}\sigma_{t(a)}(m)}
      {1\over p_{t(a)}}\bigl(
       p_{t(a)}\rk \E_{t(a)}^{(j_{t(a)})}-j_{t(a)}r_{t(a)}\bigr)
      -{j_{t(a)}\over p_{t(a)}}\\ 
    \phantom{0(\ge)\sum_{a}b_a\biggl[\q}+
    {p_{h(a)}+\sigma_{h(a)}(m)\over r_{h(a)}\sigma_{h(a)}(m)}
       {1\over p_{h(a)}}\bigl(
       p_{h(a)}\rk \E_{h(a)}^{(j_{h(a)})}-j_{h(a)}r_{h(a)}\bigr)
       +{j_{h(a)}\over p_{h(a)}}\biggr].
\end{array}
\end{equation} 
Applying Formula (\ref{weightii}) and multiplying the result by
$\sigma(m)=\sigma_1(m)\cdot...\cdot\sigma_n(m)$, (\ref{weightiii})
becomes equivalent to
\begin{equation}
\label{weightiv}
\begin{array}{ccc}
0\q(\ge)\q
\sum_{a\in A} b_a\biggl[\check{\sigma}_{t(a)}(m) \Bigl\{j_{t(a)}
-\rk\E^{(j_{t(a)})}_{t(a)}
{\chi(\E_{t(a)}(m))-\sigma_{t(a)}(m)\over \rk\E_{t(a)}}
\Bigr\}
\\
\phantom{0\q(\ge)\q\sum_{a\in A}b_a\biggl[\q}
   +\check{\sigma}_{h(a)}(m)\Bigl\{j_{h(a)}
   -\rk\E^{(j_{h(a)})}_{h(a)}
   {\chi(\E_{h(a)}(m))+\sigma_{h(a)}(m)\over \rk\E_{h(a)}}\Bigr\}
   \biggr].
\end{array}
\end{equation}
Assume first that $\catR:=(\E_i,i\in V;\phi_a,a\in A)$ is 
$\vartheta$-(semi)stable. It follows from Theorem~\ref{bound22} that
$\catR(m)$ fulfills the condition of 
$(\sigma_i(m),i\in V; b_a, a\in A)$-sectional (semi)stability for
all sub-representations $(\F_i(m),i\in V;\phi_a^\p(m), a\in A)$ for
which $\F_i(m)$ is globally generated for all $i\in V$, in particular,
for the sub-representation
$(\E_i^{(j_i)}(m),i\in V;\phi_{a|\E_{t(a)}^{(j_{t(a)})}}(m), a\in A)$.
Since $j_i\le h^0(\E_i^{(j_i)}(m))$ for all $i\in V$, this implies
(\ref{weightiv}).
\par
Now, suppose $\catg$ is (semi)stable for the given linearization
and let $(\F_i,i\in V;\phi_a^\p, a\in A)$ be a sub-representation
of $\catR$ for which $\F_i(m)$ is globally generated for all $i\in V$.
Choose bases $v_1^i,...,v_{p_i}^i$ of the $V_i$ for which there are
indices $j_i$ with $H^0(k_i(m))(V^{(j_i)}_i)=H^0(\F_i(m))$, $i\in V$.
Obviously,
$$
f_a\bigl(V_{t(a)}^{(j_{t(a)})}\bigr)\q\subset\q V_{h(a)}^{(j_{h(a)})}
\q\hbox{for every arrow } a\in A.
$$
Thus, (\ref{weightiii}) shows that 
$(\sigma_i(m),\allowbreak i\in V;b_a, a\in A)$-sectional (semi)stability is verified
for $(\F_i(m),\allowbreak i\in V;\phi_a^\p(m), a\in A)$. By Theorem~\ref{bound22}, this
implies that $\catR$ is $\vartheta$-(semi)stable.
\par
To see the assertion about the polystable points, we first remark that
a point ${\frak g}:=\mathop{\rm Gies}([k_i], i\in V; [\phi_a], a\in A)\in 
{\frak G}^{ss}$ fails to be stable if and only if there
is a destabilizing  sub-representation
$(\F_i, i\in V; \phi_a^\p, a\in A)$. As explained before, this gives
rise to a certain one parameter subgroup $\la$.
Take $\lim_{z\ra \infty}\la(z)([k_i], i\in V; [\phi_a], a\in A)$. Then,
it is not hard to see that the representation corresponding to that
point is $(\F_i,i\in V; \phi_a^\p, a\in A)
\oplus (\E_i/\F_i, i\in V; \ol{\phi}_a,a\in A)$.
From the already proven semistable version of our theorem, it follows
that this representation is again semistable and that
$\lim_{z\ra \infty}\la(z){\frak g}$ is a semistable point.
Hence, it is clear that ${\frak g}$ will be a polystable point
if and only if every destabilizing sub-representation
of $\catR$  is a direct summand, or, in other words,
$\catR$ is $\vartheta$-polystable.
\par
ii) 
We will follow \cite{HL}, Prop.~4.4.2, and 
use the valuative criterion of properness.
Let $(C,0):=\Spec R$ where $R$ is a discrete valuation ring.
By assumption, we are given a map $h\colon C\lra
{\frak G}^{ss}$
which lifts over
$C\setminus\{0\}$ to ${\frak T}^{\vartheta-ss}_{\catp}$.
In particular, by the universal property of 
${\frak T}^{\vartheta-ss}_{\catp}$, there is a
family
$$
\bigl(V_i\otimes \pi_X^*\O_X(-m)\lra{\frak
E}_i^0,i\in V;\phi^0_a, a\in A\bigr)
$$
parametrized by $C\setminus\{0\}$.
This can be extended to a certain family
$$
\bigl(V_i\otimes \pi_X^*\O_X(-m)\lra\widetilde{\frak
E}_i,i\in V;\widetilde{\phi}_a, a\in A\bigr)
$$
on $C\times X$. Here, $V_i
\otimes \pi_X^*\O_X(-m)\lra\widetilde{\frak
E}_i$ are families of not necessarily torsion free
quotients with Hilbert polynomial $P_i$, $i\in V$,
and $\widetilde{\phi}_a\in\Hom(
\widetilde{\frak
E}_{t(a)}, \widetilde{\frak E}_{h(a)})$.
For each arrow $a\in A$, there is a commutative diagram
$$
\begin{CD}
0 @>>> \widetilde{\frak E}_{t(a)} @>>> {\widetilde{\frak
E}}_{t(a)}^{\vee\vee} @>>>  {\frak N}_{t(a)} @>>> 0\phantom{.}\\
@| @V{\widetilde{\phi}}_aVV
@V\widetilde{\phi}_a^{\vee\vee}VV
@V\overline{\phi}_aVV @|\phantom{.}\\
0 @>>> \widetilde{\frak E}_{h(a)} @>>> {\widetilde{\frak
E}_{h(a)}}^{\vee\vee} @>>> {\frak N}_{h(a)} @>>> 0.
\end{CD}
$$
As in \cite{HL}, define ${\frak N}_i^\p\subset
{\frak N}_i
$, $i\in V$, as the union
of the kernels of the multiplications by $t^n$, $n\in\N$,
$t$ a generator of the maximal ideal of $R$.
Next, we set ${\frak E}_{i}:=\ker\bigl(\widetilde{\frak
E}^{\vee\vee}_{i}\lra {\frak N}_{i}/{\frak N}^\p_{i}\bigr)$, $i\in V$.
These are $C$-flat families of torsion free coherent sheaves on $C\times X$.
Since $\overline{\phi}_a$ maps  ${\frak N}^\p_{t(a)}$ to 
${\frak N}^\p_{h(a)}$,
the map $\widetilde{\phi}_a^{\vee\vee}$ induces a homomorphism
$\phi_{C,a}\colon {\frak E}_{t(a)}\lra {\frak E}_{h(a)}$.
By construction, there are homomorphisms ${\frak
k}_i\colon V_i\otimes\pi_X^*\O_X(-m)\lra
{\frak
E}_i$, $i\in V$, which coincide
on $(C\setminus\{0\})\times X$ with the quotients we started with and
which become generically surjective when restricted to $\{0\}\times X$.
The family
$$
\bigl({\frak k}_i\colon V\otimes\pi_X^*\O_X(-m)\lra {\frak E}_i,i\in V;
\phi_{C,a}, a\in A\bigr)
$$
defines a morphism of $C$ to
${\frak G}$
which is, of course, the map $h$ of the beginning.
\par
Let $k_i\colon V_i\otimes\O_X(-m)\lra \E_i$ be the restriction
of ${\frak k}_i$ to $\{0\}\times X$, $i\in V$, and $[f_a\colon
V_{t(a)}\lra V_{h(a)}, a\in A]$ the ${\frak P}$-component
of $h(0)$.
We claim that $H^0(k_i(m))$ must be injective for all
$i\in V$.
To see this, set $K_i:=\ker(V_i\lra H^0(\E_i(m)))$, $i\in V$,
and assume that not all the $K_i$ are trivial.
For each $i$, let $v^i_1,...,v_{j_i}^i$ be a basis for $K_i$ and
complete it to a basis $v^i_1,...,v_{p_i}^i$ of $V_i$.
It follows from the construction that $f_a(K_{t(a)})\subset K_{h(a)}$
for all arrows $a\in A$.
Therefore, evaluating the semistability condition
yields
\begin{eqnarray*}
0 &\ge & \sum_{a\in A} b_a\bigl(\check{\sigma}_{t(a)}(m)
j_{t(a)}+\q\check{\sigma}_{h(a)}(m)
j_{h(a)}\bigr)
\end{eqnarray*}
which is impossible.
\par
Using now $H_i:=k_i(V_i)\subset H^0(\E_i(m))$, one can check with the
same methods as
before that $(\E_i(m),i\in V; \phi_{C,a|\{0\}\times X}(m), a\in A)$
is sectional semistable w.r.t.\ the
parameters $(\sigma_i(m),i\in V; b_a, a\in A)$.
But this implies $h^0(\E_i(m))=p_i$ for $i\in V$, so that all
the ${\frak k}_i$ are honest quotients.
This means that
$$
\bigl({\frak k}_i\colon V\otimes\pi_X^*\O_X(-m)\lra {\frak E}_i,i\in V;
\phi_{C,a}, a\in A\bigr)
$$
defines a morphism $C\lra {\frak T}$ which maps by our previous
calculations to ${\frak T}^{\vartheta-ss}_{\catp}$, 
thus providing the desired lifting of $h$.

\bigskip
\par

Universit\"at GH Essen\par
FB 6 Mathematik und Informatik\par
D-45117 Essen\par
Deutschland\par
\smallskip
\it E-Mail: \tt alexander.schmitt\@@uni-essen.de
\end{document}